\documentclass[a4paper]{article}

\usepackage{booktabs}
\usepackage[english]{babel}
\usepackage[utf8]{inputenc}
\usepackage{amssymb}
\usepackage[fleqn]{amsmath}
\usepackage{arydshln}
\usepackage{amsthm}
\usepackage{graphicx}
\usepackage{framed}
\usepackage[colorinlistoftodos]{todonotes}
\usepackage{subcaption}
\usepackage{multirow}
\usepackage{enumitem}
\usepackage{titling} 
\usepackage{titlesec}
\usepackage[top=1.25in, bottom=1.25in, left=1.25in, right=1.25in]{geometry}
\usepackage{breqn}
\usepackage[numbers]{natbib}
\usepackage{tikz}
\usepackage{ulem}
\usepackage{caption}
\usetikzlibrary{shapes,snakes}
\usepackage{stackengine}
\usetikzlibrary{arrows,%
	petri,%
	topaths}%
\usepackage{algorithm}
\usepackage{algpseudocode}
\usepackage{tabstackengine}

\TABbinary
\newcounter{example}
\newenvironment{example}[2]{
\refstepcounter{example}
\label{#1}\begin{leftbar}\textbf{\noindent \hspace{-5.3mm}Example \theexample} \\#2\end{leftbar}
}{}

\newcounter{drafts}

\newcounter{module}
\newcounter{model}
\newcommand \linedabstractkw[2]{
  \renewcommand\maketitlehookd{%
    \mbox{}\medskip\par
    \centering
    \hrule\medskip
    \begin{minipage}{0.92\textwidth}
    \textbf{Abstract}\\ #1\\
    
    \textit{Keywords: }#2
    \end{minipage}\medskip\hrule\medskip
    }      
}

\newcommand \ERPWtemplate[3]{
\usepackage{fancyhdr}
\pagestyle{fancy}
\fancyhf{}
\fancyhead[RO]{#1}
\fancyhead[LO]{#2}
\fancyfoot[CO]{#3}
\fancyfoot[RO]{\thepage}

\renewcommand{\headrulewidth}{1pt}
\renewcommand{\footrulewidth}{1pt}
}

\stackMath
\newcommand{\trace}[1]{\text{Tr }#1}

\newcommand{\im}[1]{\text{Im}\left(#1\right)}
\newcommand{\re}[1]{\text{Re}\left(#1\right)}
\newcommand{\rank}[1]{\text{Rank}\left(#1\right)}
\newcommand{\rk}[1]{\text{Rank}~#1}

\newcommand{\E}{\mathcal{E}}

\newcommand{\w}{\omega}
\newcommand{\ww}{\overline{\omega}}

\newcommand{\bG}{\overline{G}}
\newcommand{\bE}{\overline{E}}

\theoremstyle{plain}
\newtheorem{theorem}{Theorem}[section]

\theoremstyle{plain}
\newtheorem{lemma}[theorem]{Lemma}

\theoremstyle{plain}
\newtheorem{proposition}[theorem]{Proposition}

\theoremstyle{remark}
\newtheorem{remark}{Remark}
\theoremstyle{question}
\newtheorem{question}{Question}

\theoremstyle{plain}
\newtheorem{corollary}{Corollary}[theorem]

\theoremstyle{definition}
\newtheorem{definition}{Definition}

\theoremstyle{plain}

\newcommand{\ch}{\mathcal{H}}

\DeclareRobustCommand{\VAN}[3]{#2} 

\ERPWtemplate{P. Wissing \& E.R. van Dam}{Signed Directed Graphs}{\leftmark}

\tikzset{group/.style = {shape=circle,draw,dotted,minimum size=1em}}
\tikzset{vertex/.style = {shape=circle,draw,minimum size=1em}}
\tikzset{arc/.style = {->,> = latex'}}
\tikzset{edge/.style = {-,> = latex'}}
\tikzset{negarc/.style = {->,> = latex',dashed}}
\tikzset{negedge/.style = {-,> = latex',dashed}}
\tikzset{tree/.style = {-,> = latex',line width=.7mm}}

\title{
Spectral Fundamentals and Characterizations of Signed Directed Graphs
}
\author{
Pepijn Wissing\thanks{Corresponding author: p.wissing@tilburguniversity.edu},~ Edwin R. van Dam \\ \small{Department of Econometrics and Operations Research, Tilburg University}
}

\begin{document}
\linedabstractkw{
The spectral properties of signed directed graphs, which may be naturally obtained by assigning a sign to each edge of a directed graph,  have received substantially less attention than those of their undirected and/or unsigned counterparts. 
To represent such signed directed graphs, we use a striking equivalence to $\mathbb{T}_6$-gain graphs to formulate a Hermitian adjacency matrix, whose entries are the unit Eisenstein integers $\exp(k\pi i/3),$ $k\in \mathbb{Z}_6.$
Many well-known results, such as (gain) switching and eigenvalue interlacing, naturally carry over to this paradigm. 
We show that non-empty signed directed graphs whose spectra occur uniquely, up to isomorphism, do not exist, but we provide several infinite families whose spectra occur uniquely up to switching equivalence. 
Intermediate results include a classification of all signed digraphs with rank $2,3$, and a deep discussion of signed digraphs with extremely few (1 or 2) non-negative (eq. non-positive) eigenvalues.
}{
Complex unit gain graphs, Hermitian, Spectra of digraphs, Signed graphs
}
\maketitle

\section{Introduction}                       
Eigenvalues of graphs have fascinated researchers in various fields for decades. 
The various representations of (undirected) graphs have lead to numerous results, connecting spectral properties of a graph to its structural characteristics. 
A natural question, that particularly intrigues the authors, is whether or not a graph may be uniquely constructed from a given spectrum.
In the standard undirected graph paradigm, many results relating to this question are known. 
For a concise survey on the topic, the interested reader is referred to \cite{vandam2003, vandam2009}. 


In recent years, various authors (see, e.g., \cite{ reff2012spectral, mehatari2018adjacency, samanta2019}) have considered the spectral properties of so-called \textit{complex unit gain graphs}.
Let $\mathbb{T} = \{z\in\mathbb{C}~:~ |z|=1\}$ be the multiplicative group of unit complex numbers.
Then, briefly put, such graphs may be considered to be weighted, symmetric directed graphs, whose weights $\psi(u,v)\in\mathbb{T}$ satisfy $\psi(u,v) = \psi(v,u)^{-1}$ for all adjacent pairs of vertices $u,v$. 
The corresponding gain matrices, which represent gain graphs in the natural way, are Hermitian and therefore satisfy a number of convenient spectral properties.  

The spectral properties of various special cases of these complex unit gain graphs have been previously studied. 
In particular, signed graphs \cite{zaslavsky1982signed,belardo2014balancedness} and the Hermitian adjacency matrices $H$ and $N$ for mixed graphs \cite{guo2017hermitian, liu2015hermitian, mohar2020new} come to mind. 
Such cases have, in hindsight, simply allowed subsets of $\mathbb{T}$ and accompanied the relevant entries with an appropriate interpretation. 
For example, in \cite{mohar2020new}, the allowed entries are $\{\exp(ik\pi/3)~|~k\in \{-1,0,1\}\}$, where the strictly real entry represents the a symmetric edge, and those with a nonzero imaginary part represent an oriented edge. 
This work generalizes that line, allowing gains in the subgroup $\{\exp(ik\pi/3)~|~k\in [6]\}=:\mathbb{T}_6.$ 
Accordingly, the corresponding matrices naturally represent \textit{signed directed graphs}.

Albeit quietly, these objects often play a central role in the modeling of dynamic systems.
In this context, a dedicated group of researchers (see, e.g., \cite{catral2009allow,hall2006sign,driessche2018sign}) has been studying signed digraphs stemming from \textit{sign-pattern matrices}; effectively imposing no assumptions on the conversion rate between states in such a system, other than the sign of the effect. 
Conclusions regarding, e.g., the stability of a (or, sometimes, any) system with the corresponding signs are obtained using spectral techniques.
In contrast, this article focuses on the above-mentioned gain graph representation of a signed directed graph, and the spectral consequences thereof.

Our ultimate interest lies with the question whether or not this representation of signed directed graphs may offer sufficient combinatorial information, in order to uniquely determine such a signed directed graph, when given its spectrum. 
The answer to this question consists of two parts.
First, we show that any signed directed graph has a partner that is switching equivalent (and thus cospectral), but not isomorphic. 
This is a natural consequence of the fact that $\mathbb{T}_6$ is closed under multiplication, due to which it is always possible to apply gain switching and obtain a non-isomorphic signed directed graph. 
We do, however, obtain several families of signed directed graphs that are switching equivalent to every signed directed graph to which they are cospectral. 

In order to streamline the discussion, we first set out to classify signed directed graphs that satisfy particular spectral conditions.
We show that one may without loss of generality assume the gains of a spanning tree, which aids the classification significantly.
By applying this idea, in conjunction with eigenvalue interlacing we classify all signed directed graphs whose rank is $2$ or $3$; we find a notably concise characterization of all such signed digraphs, which may be described as twin expansions of either an edge, a triangle, or the transitive tournament of order four.

Subsequently, we provide an extensive discussion of signed directed graphs with exactly $1$ or $2$ non-negative eigenvalues. 
We show that such graphs are highly dense and provide a list of necessary properties, though already in the case of $2$ non-negative eigenvalues, the complete list of candidates quickly becomes unwieldy. 
Thus, we focus on a few special cases, such as clique expansions of the $5$-cycle and the $4$-path. 
In particular, we characterize all signed directed graphs on these minimally dense graphs, such that the resulting signed digraphs admit to the imposed requirements.

The above characterizations are then used to consider spectral determination. 
Through a series of counterexamples, we show that the discussed low rank signed digraphs are not, in general, determined by their spectra. 
However, by applying a sequence of counting arguments to the lists obtained above, we are then able to prove that, among others, several of the families with $2$ non-negative eigenvalues are determined by their spectrum. 
Specifically, in addition to a number of sporadic examples, we find several arbitrarily large graphs, obtained as clique expansions of $C_4, P_4$ or $C_5$, that admit signed digraphs cospectral only to switching equivalent signed digraphs. 

The contents of this paper are organized as follows. 
We first provide a thorough introduction of the subject matter, in Section \ref{sec: prelim}.  
In Section \ref{sec: elementary}, we showcase adaptations of well-known results from the adjacent fields, that more or less carry over to the current framework. 
Sections \ref{sec: low rank} and \ref{sec: extreme ev} are concerned with the characterization of signed digraphs that satisfy an imposed set of spectral requirements.  
The obtained knowledge is then applied in Section \ref{sec: cospectrality}, to investigate spectral characterizations of signed digraphs.  
Finally, we conclude with a collection of open questions. 

The main tools used throughout are eigenvalue interlacing and expansion of graphs via lexicographic products with either empty or complete graphs.
Additionally, Lemmas \ref{lemma: edges} and \ref{lemma: triangles}, which count respectively  edges and triangles, and Proposition \ref{prop: we can always make a tree equal}, which allows us to fix the signature of a subset of the  edges without loss of generality, are frequently applied throughout.  

\section{Preliminaries}                      
\label{sec: prelim}
Let us first thoroughly define the key concepts and notation that is used throughout this work.

\subsection{Basic definitions}
A \textit{directed graph}, usually abbreviated \textit{digraph}, is a pair $D=(V,E),$ where $V$ (sometimes $V(D)$, for clarity) denotes the vertex set, whose order is $n$. Further, $E\subseteq V\times V$ denotes the  edge set, whose members are ordered pairs of vertices, called \textit{arcs} and denoted $uv$, for vertices $u,v\in V$. If $uv\in E$ and $vu\in E$, then $(u,v)$ is also called a \textit{digon.} 
To circumvent some of the more tedious technicalities, we will often use the word \textit{edge} in (signed) digraph context to indicate an edge in the underlying graph. 
Throughout, the term \textit{graph} will be reserved for undirected graphs, which will commonly be denoted with $G$.   

The objects studied in this work are, in essence, $\mathbb{T}_6$-gain graphs. 
These are \textit{complex unit gain graphs} \cite{reff2012spectral} whose gain groups are restricted to the multiplicative group $\mathbb{T}_6=\{\exp(ik\pi/3)~|~k\in[6]\}=\{\w^k~|~k\in[6]\},$ where $\w=(1+i\sqrt{3})/2$ as is denoted throughout. 
Explicitly, such a $\mathbb{T}_6$-gain graph $\Phi$ is said to be the pair $(D,\varphi),$ where $D=(V,E)$ is a symmetric\footnote{A directed graph is called symmetric when $uv$ is an arc if and only if $vu$ is an arc.} directed graph, and $\varphi: E\mapsto \mathbb{T}_6$ is a \textit{gain function} such that $\varphi(uv)=\varphi(vu)^{-1}$. 
We will often refer to $\varphi$ as the \textit{signature} of $\Phi$.

The main discussion in this paper is concerned with the matrices that are associated with the described structures.
The following matrix arises naturally. 
\begin{definition}\label{def: eisenstein matrix}
Let $\Phi = (D,\varphi)$ be a $\mathbb{T}_6$-gain graph.
Then its \textit{unit Eisenstein matrix} (or Eisenstein matrix for short) $\E$ is defined as
\[\E_{uv}(\Phi) = \begin{cases}\varphi(u,v) &\text{~if~} uv\in E\\ 0 & \text{~otherwise.}\end{cases}\]
\end{definition}
Following the definition above, the nonzero entries of $\E(\Phi)$ are exactly the unit elements of the imaginary quadratic ring $\mathbb{Z}[\omega].$
The elements of the latter group are called the Eisenstein integers \cite{greaves2012cyclotomic}, which justifies the terminology.

As noted before, an Eisenstein matrix whose entries have non-negative real parts coincides exactly with the alternative Hermitian adjacency matrix for directed graphs, proposed by \citet{mohar2020new}. 
Since the negative counterpart to every such (non-zero) entry is also contained in $\mathbb{T}_6$, an arbitrary Eisenstein matrix naturally represents a \textit{signed directed graph}.
That is, any $\mathbb{T}_6$-gain graph coincides with the natural Hermitian adjacency matrix\footnote{Note that the original Hermitian adjacency matrix  \cite{guo2017hermitian, liu2015hermitian} does not allow for a natural inclusion of the signs, whereas the variant \cite{mohar2020new} does.} of a directed graph whose edges are accompanied with a weight that is either $1$ or $-1$, and vice versa. 
Hereafter, we will use the latter perspective in our discussion, though the implications and applicable theory of the gain graph equivalent are widely applied. 

\begin{figure}[t]
\begin{center}
\begin{tikzpicture}[scale=1.3,cap=round,>=latex]
        \draw[->] (-1.5cm,0cm) -- (1.5cm,0cm) node[right,fill=white] {Re};
        \draw[->] (0cm,-1.5cm) -- (0cm,1.5cm) node[above,fill=white] {Im};

        \draw[thick] (0cm,0cm) circle(1cm);

        \foreach \x in {0,60,...,360} {
                \draw[gray] (0cm,0cm) -- (\x:1cm);
                \filldraw[black] (\x:1cm) circle(0.6pt);
        }

        \foreach \x/\xtext in {
            60/\omega^1,
            120/\omega^2,
            180/\omega^3,
            240/\omega^4,
            300/\omega^5,
            360/\omega^6}
                \draw (\x:0.65cm) node[fill=white] {$\xtext$};

        \foreach \x/\xtext in {
            60/\omega,
            120/-\bar{\omega},
            240/-\omega,
            300/\bar{\omega}}
                \draw (\x:1.3cm) node[fill=white] {$\xtext$};
                
        \draw (-1.25cm,0cm) node[above=1pt] {$-1$}
              (1.20cm,0cm)  node[above=1pt] {$1$};
    \end{tikzpicture}
    \caption{The possible entries of $\E$. 
    }
    \end{center}
    \end{figure}
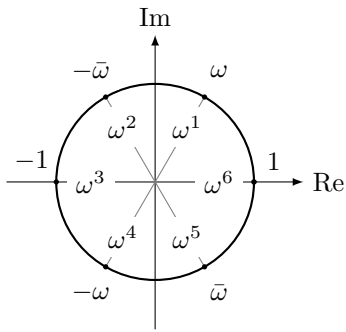

For a given signed digraph $\Phi$ with Eisenstein matrix $\E$, its \textit{characteristic polynomial} $\chi(\lambda)$, is said to be the characteristic polynomial of its Eisenstein matrix. That is, $\chi(\lambda) = \det\left(\lambda I-\E\right).$ 
The \textit{eigenvalues} of $\Phi$ are the roots of $\chi(\lambda)$; the collection of eigenvalues $\lambda_1\geq\lambda_2\geq\ldots\geq\lambda_n$ is called the \textit{spectrum}.
Two signed digraphs are said to be \textit{cospectral} if their spectra  (eq. characteristic polynomials) coincide. 

Throughout, we will encounter several arguments rooted in the concept known as eigenvalue interlacing. We include the formal definition, below.
\begin{lemma}\label{lemma: interlacing} \citep{haemers1995, godsil2001algebraic}
Suppose $A$ is a Hermitian $n\times n$ matrix with eigenvalues $\lambda_1\geq \ldots\geq \lambda_n$. Then the eigenvalues $\mu_1\geq \ldots\geq \mu_m$ of an $m\times m$ principal submatrix of $A$ satisfy $\lambda_i \geq \mu_i \geq \lambda_{n-m+i}$ for $i\in[m]$.
\end{lemma}
We recall some common notions. 
For a given signed digraph $\Phi$, its underlying graph is obtained by replacing every arc with a digon, and setting all signs to $1$. 
We define the \textit{underlying graph operator} $\Gamma(\cdot)$, that maps a signed digraph to its underlying graph, and we use $\cup$ to denote the disjoint union of two (signed) (di)graphs $A$ and $B$ as $A\cup B$.
For a given (signed) (di)graph $\Phi$, we denote the subgraph induced by $U\subseteq V(\Phi)$ as $\Phi[U]$.
Two (signed) (di)graphs $A$ and $B$ are said to be \textit{isomorphic} (denoted $A\cong B$) if $\E(A)=P\E(B)P^\top$, for a permutation matrix $P$; that is, if they are equal up to a relabeling of the vertices.
We will often denote $[k]=\{1,\ldots,k\}$.

The gain of a walk is the product of the gain of the traversed edges. 
If a walk is circular and such that the only repeated vertices are the first and the last, it is said to be a cycle, commonly denoted $C$. 
Of particular interest for the presented spectral analysis is the gain of a cycle, denoted $\varphi(C)$. 
Clearly, $\varphi(C)\in\mathbb{T}_6$ for any $C$.
A cycle $C$ is called \textit{real} if $\im{\varphi(C)}=0$, and it is called \textit{positive} (resp. \textit{negative}) if $\re{\varphi(C)}>0$ (resp. $<0$). 
It should be noted that that the gain of a cycle is, strictly speaking, dependent on the direction in which it is traversed. 
Specifically, traversing the cycle in the converse direction yields the complex conjugate of the gain value. 
However, since the real part contains all of the information that is interesting from a spectral point of view (see Theorem \ref{thm: coefficients}), the choice of traversal direction is, for the purposes of this work, inconsequential and thus not specifically mentioned.

To end this section, we include a list of graphs and the names by which they are indicated.
As usual, a complete graph is denoted $K_n$, and a complete $k$-partite graph is denoted $K_{n_1,\ldots,n_k}$.
Further, the empty graph  is denoted $O_n$, the path is denoted $P_n$, and a cycle is denoted $C_n$. 
The \textit{transitive tournament} $T_n$ is the digraph whose  arc set is exactly $E(T_n) = \left\{uv ~|~ u\leq v\text{~for~}u,v\in [n]\right\}.$
The remaining two graphs that are occasionally referred to are illustrated in Figure \ref{fig: small graphs}. 
Lastly, a brief summary of drawing conventions is included in Table \ref{tab: drawing conventions}.

\begin{table}[b]
\begin{minipage}{0.51\textwidth}\centering
    \begin{minipage}{.45\textwidth}\centering
     \begin{tikzpicture} 
     \node[vertex] (1) at (0,0) {};
     \node[vertex] (2) at (0,1) {};
     \node[vertex] (3) at (1,0.5) {};
     \node[vertex] (4) at (2,0.5) {};
     
     \draw[edge] (1) to node{} (2);
     \draw[edge] (1) to node{} (3);
     \draw[edge] (3) to node{} (2);
     \draw[edge] (3) to node{} (4);
     
     \end{tikzpicture}\\
     (a) 3-pan
     \end{minipage}
     \begin{minipage}{.45\textwidth}\centering
     \begin{tikzpicture} 
     
     \node[vertex] (1) at (0,0) {};
     \node[vertex] (2) at (.75,0.7) {};
     \node[vertex] (3) at (0.3,1) {};
     \node[vertex] (4) at (-0.3,1) {};
     \node[vertex] (5) at (-.75,0.7) {};

     \draw[edge] (1) to node{} (2);
     \draw[edge] (1) to node{} (3);
     \draw[edge] (1) to node{} (4);
     \draw[edge] (1) to node{} (5);
     \draw[edge] (2) to node{} (3);
     \draw[edge] (3) to node{} (4);
     \draw[edge] (4) to node{} (5);
     
     \end{tikzpicture}\\
     (b) Gem
     \end{minipage}

     \captionof{figure}{Two small graphs and their names.}
     \label{fig: small graphs}
\end{minipage}
\begin{minipage}{.48\textwidth}
\centering
\begin{tabular}{llr}
\toprule
Edge type     & Drawing & $\E_{lr}$ \\ \midrule
Positive digon &       \begin{tikzpicture}
        \node[vertex] (1) at (0,0) {};
        \node[vertex] (2) at (1,0) {};
        \draw[edge] (1) to node{} (2);
    \end{tikzpicture}  & $1$    \\
Negative digon &   \begin{tikzpicture}
        \node[vertex] (1) at (0,0) {};
        \node[vertex] (2) at (1,0) {};
        \draw[negedge] (1) to node{} (2);
    \end{tikzpicture}   & $-1$     \\
Positive arc  &    \begin{tikzpicture}
        \node[vertex] (1) at (0,0) {};
        \node[vertex] (2) at (1,0) {};
        \draw[arc] (1) to node{} (2);
    \end{tikzpicture}    & $\w$   \\
Negative arc  &    \begin{tikzpicture}
        \node[vertex] (1) at (0,0) {};
        \node[vertex] (2) at (1,0) {};
        \draw[negarc] (1) to node{} (2);
    \end{tikzpicture}   & $-\w$    \\ \bottomrule
\end{tabular}
\captionof{table}{Drawing conventions. Here, $l$ and $r$ are the left and right vertices, respectively.}
\label{tab: drawing conventions}
\end{minipage}
\end{table}
\subsection{Expansions}
In Sections \ref{sec: low rank} and \ref{sec: extreme ev}, we will be looking to construct arbitrarily large signed digraphs, based on smaller structures that we know admit to some predetermined set of requirements. 
Depending on the context, we will be looking to add either \textit{twins} or \textit{pseudotwins} to a signed digraph. 
While these conceptual ideas are widely known, we include a formal definition, as the details tend to vary. 
In particular, we include gain switching, details about which will follow shortly, into these definitions to forego some tedious discussion. 
\begin{definition}
\label{def: twins}
Let $\Phi=(G,\varphi)$ be a signed digraph of order $n$, whose Eisenstein matrix is $\E,$ and let $u,v\in V(\Phi)$ be distinct vertices. 
If, for some gain switching matrix\footnote{The formal definition is provided in Section \ref{sec: switching equivalence}.} $X$ and all $z\in V$ we have $\E_{uz}=(X\E X^{-1})_{vz}$, then $u$ and $v$ are called \textit{(switching) twins}.
\end{definition}
\begin{definition}
\label{def: pseudotwins}
Let $\Phi=(G,\varphi)$ be a signed digraph of order $n$, whose Eisenstein matrix is $\E,$ and let $u,v\in V(\Phi)$ be distinct vertices. 
If, for some gain switching matrix $X$ and all $z\in V$ we have $(\E+I)_{uz}=(X\E X^{-1} + I)_{vz}$, then $u$ and $v$ are called \textit{(switching) pseudotwins}.
\end{definition}

In essence, two nodes $u$ and $v$ are twins or pseudotwins if their respective relations to the remaining vertices in $V$ are equivalent.
The former additionally requires $u$ and $v$ to be non-adjacent, while the latter requires that they are adjacent and additionally that all triangles containing both $u$ and $v$ have gain $1$. 

Now, we may straightforwardly define expansion and reduction operators that respectively grow and shrink signed digraphs, while preserving the underlying structure. 
\begin{definition}\label{def: twin expansion}
Let $\Phi=(G,\varphi)$ be a signed digraph of order $n$, and let $\tau\in\mathbb{N}^n$ be a vector. Then the \textit{twin expansion of $\Phi$ with respect to $\tau$}, denoted $TE(\Phi,\tau)$, is obtained by introducing $\tau_j$ twins of vertex $j$ into $\Phi$. Conversely, the \textit{twin reduction} $TR(\Phi)$ of $\Phi$ is obtained by removing all but one of every collection of twins from $\Phi$. Accordingly, $\Phi$ is called \textit{twin reduced} if $\Phi=TR(\Phi)$. 
\end{definition}
\begin{definition}\label{def: clique expansion}
Let $\Phi=(G,\varphi)$ be a signed digraph of order $n$, and let $\tau\in\mathbb{N}^n$ be a vector. Then the \textit{clique expansion of $\Phi$ with respect to $\tau$}, denoted $CE(\Phi,\tau)$, is obtained by introducing $\tau_j$ pseudotwins of vertex $j$ into $\Phi$. The \textit{clique reduction} $CR(\Phi)$ of $\Phi$ is obtained by removing all but one of every collection of pseudotwins from $\Phi$. Accordingly, $\Phi$ is called \textit{clique reduced} if $\Phi=CR(\Phi)$. 
\end{definition}

Alternatively, one may think of the expansion operations defined above as taking the lexicographic product (see, e.g., \cite{godsil2001algebraic}) of a signed digraph $\Phi$ with the collection $\{O_{\tau_1},O_{\tau_2},\ldots,O_{\tau_n}\}$, in the case of twin expansion, or $\{K_{\tau_1},K_{\tau_1},\ldots,K_{\tau_1}\}$ in the case of clique expansion. 
Note that the number of nonzero (resp. not $-1$) eigenvalues is unaffected by the above expansion operators. 
Other authors (e.g., \cite{mohar2016Hermitian, oboudi2016characterization}) have used the concepts above with varying notation; the authors prefer the definition in terms of operators, to make the distinction between them clear. 

\begin{remark}
Since, for both expansion operators, vertex $j$ is mapped to a group of $\tau_j$ vertices, the ordering of $\tau$ and the corresponding labeling of the graph that is to be expanded both matter.
Without explicit mention hereafter, we will always label the vertices of a path graph such that its edges are $(1,2),(2,3),\ldots,(n-1,n),$ and the vertices of a cycle graph such that its edges are $(1,2),(2,3),\ldots,(n-1,n),(n,1).$ Other cases will be explicitly illustrated.
\end{remark}
\section{Elementary results}                 
\label{sec: elementary}
The first natural question to ask is, of course, which of the known results from related fields carry over and in what capacity. 
In this section, we will discuss a number of them, that are especially relevant for the remainder of this work. 
Additionally, we show that one may without loss of generality chose the edge gains of a spanning tree, and briefly discuss signed digraphs with symmetric spectra. 

\subsection{Counting substructures}
\label{sec: counting substructures}
A well-known result in spectral graph theory is that the number of closed walks in a graph of a given length are, in a sense, counted by the sum of its eigenvalues, exponentiated to the corresponding power. 
With respect to $\E$, a direct analogue of this idea holds. 
\begin{lemma}\label{lemma: edges}
If $\Phi$ is a signed digraph such that $\Gamma(\Phi)$ contains $m$ edges. Then
$\trace{\E(\Phi)^2} = 2m$.
\end{lemma}
\begin{proof}
Let $\bf{e}_j$ denote the columns of $\E:=\E(\Phi).$ Then $\E_{jj}^2 = \mathbf{e}_j^*\mathbf{e}_j 
= d_j,$ where $d_j$ is the degree of node $j$ in $\Gamma(\Phi).$ Hence, we have $\trace \E^2 = \sum_{j=1}^n d_j = 2m.$
\end{proof}
Note that we may categorize the three-cycles into four categories, based on the real part of their gains. 
The following then straightforwardly follows.
\begin{lemma}\label{lemma: triangles}Let $\Phi$ be a signed digraph, and let $s_{(z)}$ denote the number of triangles $t$ with $\re{\varphi(t)}=z$ that are contained in $\Phi$ as induced subdigraphs. Then \[\trace \E(\Phi)^3 = 6s_{(1)} + 3s_{(1/2)} - 3s_{(-1/2)} - 6s_{(-1)}.\]
\end{lemma}
\begin{proof}
Let $u\in V(\Phi)$ and let $\Delta_u$ be the collection of triangles in $\Gamma(\Phi)$ that contain $u$. Then, we have
\[(\E^3)_{uu} = \sum_{t\in \Delta_u}\varphi(t) = \sum_{z\in\left\{\pm 1, \pm \frac{1}{2}\right\}} 2z\cdot c_{(z)},
\]
where $c_{(z)}$ denotes the number of triangles $t$ with $\text{Re }\varphi(t)=z$ in $\Phi$ contain $u$. 
Here, the second equality holds since every triangle is traversed in two directions.
Specifically, recall that the gains of such mirror traversals are each others complex conjugate, and $\alpha + \bar{\alpha} = 2\text{Re }\alpha$ for $\alpha\in\mathbb{C}$. 
The claim then follows, since every triangle counted thrice: once for every vertex it contains. 
\end{proof}

\subsection{Coefficients theorem}
A well-known tool in spectral graph theory is known as the Harary-Sachs coefficients theorem \cite{cvetkovic1980spectra}, which computes the coefficients of the characteristic polynomial of a graph based on some structural properties of its elementary spanning subgraphs. 
For the current work, we may specialize a result by \citet{samanta2019}, to be significantly more combinatorially approachable. 
We first impose some simple notation.

A graph $G$ is called an elementary graph if each of its components is either an edge or a cycle. 
Let $\ch_j(G)$ denote the collection of all elementary subgraphs of a graph $G$ with $j$ vertices. 
For any $H\in\ch_j(G)$, let $\mathcal{C}_H$ denote the collection of all cycles in $H$. 
Let $n(H)$ and $z(H)$ respectively denote the number of negative cycles and non-real cycles in $\mathcal{C}_H$.
Finally, let $p(H)$ and $c(H)$ be the number of components and the number of cycles in $H$, respectively. 
Then, we may write the following version of said theorem.
\begin{theorem}\label{thm: coefficients}
Let $\Phi$ be a signed digraph with underlying graph $G$. Let $\chi(\lambda) = \lambda^n + a_1\lambda^{n-1} + \ldots + a_n$ be the characteristic polynomial of $\Phi$. Then the coefficients $a_j$ may be calculated as
\[
    a_j = \sum_{H\in\ch_j(G)} (-1)^{p(H) + n(H)}2^{c(H)-z(H)}.
\]
\end{theorem}
\begin{proof}
The result follows as an easy adaptation from \citet[Thm 2.7]{samanta2019}, using the following two observations: (i) any cycle $C$ has $|\re{\varphi(C)}|= 1$ if $C$ is real and $\frac{1}{2}$ otherwise, and (ii) we have (with slight abuse of notation)
\[\text{sign}\left(\prod_{C\in\mathcal{C}_H}\re{\varphi(C)}\right)= \prod_{C\in\mathcal{C}_H}\text{sign}(C) =(-1)^{n(H)}
\implies 
\prod_{C\in\mathcal{C}_H}\re{\varphi(C)} = (-1)^{n(H)} 2^{-z(H)}\]
and the result follows. 
\end{proof}


%
%

\subsection{Switching equivalence and isomorphism}
\label{sec: switching equivalence}

As we aim to classify signed digraphs whose spectra determine them, it is natural to weigh which 'other' signed digraphs are sufficiently closely related to be considered equivalent. 
In line with the literature \cite{zaslavsky1989biased}, we consider two signed digraphs to be equivalent if they are \textit{switching isomorphic.} 
First, we recall the essential switching operation\footnote{Originally due to \citet{zaslavsky1989biased}, this is essentially not a matrix-based operation. The current formulation is kept, in line with e.g. \citet{guo2017hermitian}.}.
\begin{definition}
\label{def: switching eq}
Let $\Phi$ and $\Phi'$ be signed digraphs of order $n$, whose Eisenstein matrices are $\E$ and $\E'$, respectively. 
We say that $\Phi'$ is obtained from $\Phi$ by a \textit{gain switching} if
\begin{equation}\E'= X\E X^{-1},\label{eq: switching eq}\end{equation} 
where $X$ is a diagonal matrix with $X_{ii}\in\mathbb{T}_6$ for all $i\in[n]$.
\end{definition}
By appending the above switching with vertex relabeling and taking the converse (both are well-known to yield somewhat trivial cospectral mates), the following\footnote{While functionally analogous to 'switching equivalence' \citep{guo2017hermitian,  reff2016oriented, samanta2019}, the inclusion of isomorphism is natural and more complete.} is a natural notion of equivalence.
\begin{definition}
Two signed digraphs $\Phi_1$ and $\Phi_2$ are said to be
\textit{switching isomorphic} (denoted $\Phi_1\sim\Phi_2$) if one may be obtained from the other by a sequence of diagonal switches and vertex permutations, possibly followed by taking the converse. 
\end{definition}

A well-known fact concerning switching isomorphism is the following. 

\begin{lemma}\label{lemma: sq iso => cospectral}
Any two switching isomorphic signed digraphs are cospectral.
\end{lemma}
As usual, the reverse need not be true. The saltire pair, while originally used in unsigned graph context, serves as an example for this conclusion.

It should be clear that any two switching isomorphic signed digraphs share a common underlying graph. 
Again, the reverse is in general not true. 
Consider, for instance, $\Phi = (K_n, +)$ and $-\Phi = (K_n,-),$ for $n\geq 3$. 
While their underlying graphs are clearly equal, their respective spectra are distinct and, by Lemma \ref{lemma: sq iso => cospectral}, this pair is inherently not switching isomorphic.

Naturally, we would like to have a way to conclusively determine whether or not a given pair of signed digraphs is switching isomorphic, which is found in \citep{reff2016oriented}.
Specifically, \citet{reff2016oriented} presents a sufficient condition, based on the idea (originally from \cite{zaslavsky1989biased}) that the gain of a given cycle in a gain graph is not affected by switching operations. 
This condition is, in fact, also necessary; a fact previously discussed by  \citet{samanta2019}.
The authors feel that the result follows more directly than the discussion in \cite{samanta2019} suggests. 
We present the alternative proof below.

\begin{proposition}\label{prop reff}
Let $D$ and $\Phi$ be signed digraphs on the graph $G$. Then $D\sim\Phi$ if and only if there is a $D'$ isomorphic to $D$ with  $\varphi(D'[C])=\varphi(\Phi[C])$ for every cycle $C$ in $G$. 
\end{proposition}
\begin{proof}
Sufficiency is shown in \citet{reff2016oriented}, so we will only discuss necessity.
Let $\E(\Phi) = X\E(D)X^{-1},$ and let $C=\{u_1u_2,u_2u_3\ldots,u_ku_1\}$ be a cycle in $G$. Then:
\begin{align*}
    \varphi_{\Phi}(C) &= \E(\Phi)_{u_1,u_2}\E(\Phi)_{u_2,u_3}\cdots\E(\Phi)_{u_k,u_1}\\
    &= X_{u_1,u_1}\E(D)_{u_1,u_2}X_{u_2,u_2}^{-1}\cdot 
                     X_{u_2,u_2}\E(D)_{u_2,u_3}X_{u_3,u_3}^{-1} \cdots 
                     X_{u_k,u_k}\E(D)_{u_k,u_1}X_{u_1,u_1}^{-1} \\
                   &= X_{u_1,u_1}\E(D)_{u_1,u_2}\E(D)_{u_2,u_3}\cdots\E(D)_{u_k,u_1}X_{u_1,u_1}^{-1}\\
                   &= \varphi_D(C).
\end{align*}
Finally, observe that a relabeling of the vertices changes nothing except for the indices in the above, and the claim follows.
\end{proof}
Marginally expanding on the above, \citet{samanta2019} show that one only needs the fundamental cycles of a gain graph to have equal (real parts of) gains. 
Indeed, it stand to reason that if a basis of the cycle space has equal gains, then all cycles in the cycles space agree. 
To see this, one simply needs to observe that any cycle $C$ may be obtained as the symmetric difference of fundamental cycles, and note that the gain of $C$ is simply the product of the gains of these fundamental cycles. 
In this last respect, one should exercise some care, as the traversal direction does matter here, and should be chosen such that the  edges on which the fundamental cycles intersect are traversed in opposite directions.

We conclude this section with an application of the above, which shows that cospectral signed digraphs on the same underlying graph may belong to distinct switching isomorphism classes.
\example{ex: cospectral but not iso}{Consider the signed digraphs $\Phi_a$ and $\Phi_b$ in Figure \ref{fig: counterexample cospec implies sw iso}.
A quick computation of the characteristic polynomial of their respective Eisenstein matrices yields that
\[\chi_{\Phi_a}(\lambda)= \chi_{\Phi_b}(\lambda) = \lambda^6 - 8\lambda^4 + 13\lambda^2 -5. \]
However, the gains of their fundamental cycles, shown in Table \ref{table: cycle gains counterexample}, do not coincide.
Thus, $\Phi_a$ and $\Phi_b$ do not belong to the same switching isomorphism class.
As a final note, we remark that this conclusion could also have been drawn by computing the gain of the (sole) $6$-cycle in $\Phi_a$ and $\Phi_b$, though an application of Proposition \ref{prop reff} seems appropriate.
}

\begin{table}[t]
\begin{minipage}{0.5\textwidth}\centering
    \begin{minipage}{.48\textwidth}\centering
     \begin{tikzpicture}[scale=.9] 
     \def\n{6}
     \def\rad{1.3}
     \foreach \a in {1,2,...,\n}{
     \node[vertex] (\a) at (\a*360/\n: \rad cm) {\a};
     }
     \draw[edge] (1) to node{} (2);
     \draw[edge] (1) to node{} (5);
     \draw[edge] (1) to node{} (6);
     \draw[negedge] (2) to node{} (3);
     \draw[edge] (2) to node{} (4);
     \draw[edge] (3) to node{} (4);
     \draw[arc] (5) to node{} (4);
     \draw[edge] (5) to node{} (6);
     
     \end{tikzpicture}\\
     (a) $\Phi_a$
     \end{minipage}
     \begin{minipage}{.48\textwidth}\centering
     \begin{tikzpicture} [scale=.9] 
     
    
     \def\n{6}
     \def\rad{1.3}
     \foreach \a in {1,2,...,\n}{
     \node[vertex] (\a) at (\a*360/\n: \rad cm) {\a};
     }
     \draw[edge] (1) to node{} (2);
     \draw[edge] (1) to node{} (5);
     \draw[arc] (1) to node{} (6);
     \draw[edge] (2) to node{} (4);
     \draw[negarc] (3) to node{} (2);
     \draw[edge] (3) to node{} (4);
     \draw[arc] (5) to node{} (4);
     \draw[edge] (5) to node{} (6);
     
     
     \end{tikzpicture}\\
     (b) $\Phi_b$
     \end{minipage}

     \captionof{figure}{Cospectral but not switching isomorphic signed digraphs on one underlying graph.}
     \label{fig: counterexample cospec implies sw iso}
\end{minipage}\quad
\begin{minipage}{.48\textwidth}
\centering
\begin{tabular}{@{}ccc@{}}
\toprule
$U$           & $\re{\phi(\Phi_a[U])}$ & $\re{\phi(\Phi_b[U])}$ \\ \midrule
$\{2,3,4\}$   & $-1$                   & $-1/2$                 \\
$\{1,5,6\}$   & $1$                    & $1/2$                  \\
$\{1,2,4,5\}$ & $1/2$                  & $1/2$                  \\ \bottomrule
\end{tabular}\vskip9mm
\captionof{table}{Fundamental cycle gains in Figure \ref{fig: counterexample cospec implies sw iso}}
\label{table: cycle gains counterexample}
\end{minipage}
\end{table}
\subsection{Limiting degrees of freedom}
As a direct consequence of the equivalence relations discussed before, any exercise in classification of signed digraphs would encounter an abundance of seemingly distinct digraphs, that turn out to be equivalent upon closer inspection. 
Thus, it is desirable to consider ways to limit the number of possibilities that have to be considered. 
It seems particularly practical to be able to fix a subset of the  edges to a certain type, while maintaining the certainty that all switching isomorphic classes were considered. 

In the below, we will show that any switching isomorphism class on a graph $G$ contains at least one member whose  edge-induced subdigraph coincides with a fixed spanning tree $T$ of $G$. 
This idea is quite natural from a gain graph perspective, using the well-known result that appears here as Corollary \ref{cor: forest}. 
In the interest of completeness, we include a brief proof to a result that will frequently be applied, later on. 

\begin{proposition}
\label{prop: we can always make a tree equal}
Let $G$ be a graph and let $\Phi_1,\Phi_2$ be distinct signed digraphs on $G$. Let $T\subseteq E(G)$ be a spanning tree of $G$. Then, there exists a switching matrix $Y$ such that $\Phi_2',$ obtained from $\Phi_2$ as $\E(\Phi_2') = Y\E(\Phi_2)Y^{-1}$, satisfies $\Phi_2^'[T]=\Phi_1[T]$.
\end{proposition}
\begin{proof}
Consider the edge $(u,v)\in T.$ 
Since $T$ is a spanning tree, $T\setminus (u,v)$ induces an (edge-induced) subdigraph on $G$ that consists of two disjoint components, say $V_1$ (that contains $u$) and $V_2$ (that contains $v$); see Figure \ref{fig: proof make tree equal}. 
Now, let $\E(\Phi_j)_{uv}$ denote the $(u,v)$ entry of the Eisenstein matrix corresponding to $\Phi_j$, and construct the diagonal matrix $X^{(uv)}$ as
\begin{equation}
    \label{eq: construct diag X}
    X^{(uv)}_{jj}  = \begin{cases}\E(\Phi_1)_{uv}/\E(\Phi_2)_{uv} & \text{~if~} j\in V_1 \\ 1 & \text{~if~} j\in V_2\end{cases}
\end{equation}
Now, consider the switched digraph $\Phi_2'$, whose Eisenstein matrix is $\E':=X^{(uv)}\E(\Phi_2)\left(X^{(uv)}\right)^{-1}  .$
Firstly, observe that we have $\E'_{uv} = \E(\Phi_1)_{uv}$, by construction. 
Moreover, since for any $(p,q)\in T\setminus (u,v)$ it holds that either $\{p,q\}\subset V_1$ or $\{p,q\}\subset V_2$, it follows that $\E'_{pq}=\E(\Phi_2)_{pq}$.
In other words, the arcs in $\Phi$ corresponding to exactly one edge in $T$, namely $(u,v),$ were changed by the switching with $X^{(uv)}.$
It follows that  $Y := \prod_{(u,v)\in T} X^{(uv)}$ satisfies the desired requirements. 
\end{proof}

\begin{table}[t]
\begin{minipage}{.43\textwidth}
\centering
\begin{tikzpicture}
            \node[vertex] (0) at (0,0) {$u$};
            \node[vertex] (1) at (1,1) {$v$};
            \node[vertex] (2) at (1,0) {};
            \node[vertex] (3) at (2,1) {};
            \node[vertex] (4) at (2,0) {};
            
            \draw[tree] (0) to node{} (1);
            \draw[tree] (1) to node{} (3);
            \draw[tree] (0) to node{} (2);
            \draw[tree] (2) to node{} (4);
            \draw[edge] (1) to node{} (2);
            \draw[edge] (3) to node{} (4);
            
            \draw[dotted] (1.5,1) ellipse (1.1 and .4);
            \draw[dotted] (1,0) ellipse (1.7 and .5);
            
            \draw (-.5,.5) node {$V_1$};
            \draw (2.8,1) node {$V_2$};
    \end{tikzpicture}
    
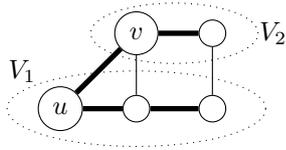
\captionof{figure}{An example graph $G$ for Prop. \ref{prop: we can always make a tree equal}. Here, the thick lines represent $T$.}
     \label{fig: proof make tree equal}
\end{minipage}\quad\quad
\begin{minipage}{0.50\textwidth}\centering
    \begin{minipage}{.48\textwidth}\centering
     \begin{tikzpicture} 
     \node[vertex] (1) at (0,1) {};
     \node[vertex] (2) at (1,1) {};
     \node[vertex] (3) at (2,1) {};
     \node[vertex] (4) at (2,0) {};
     \node[vertex] (5) at (1,0) {};
     \node[vertex] (6) at (0,0) {};
     
     \draw[tree] (1) to node{} (2);
     \draw[tree] (2) to node{} (3);
     \draw[tree] (3) to node{} (4);
     \draw[tree] (4) to node{} (5);
     \draw[tree] (5) to node{} (6);
     \draw[edge] (2) to node{} (5);
     \draw[negedge] (1) to node{} (6);
     
     \end{tikzpicture}\\
     (a) $\Phi$
     \end{minipage}
     \begin{minipage}{.48\textwidth}\centering
     \begin{tikzpicture} 
     
    \node[vertex] (1) at (0,1) {};
     \node[vertex] (2) at (1,1) {};
     \node[vertex] (3) at (2,1) {};
     \node[vertex] (4) at (2,0) {};
     \node[vertex] (5) at (1,0) {};
     \node[vertex] (6) at (0,0) {};
     
     \draw[tree] (1) to node{} (2);
     \draw[tree] (2) to node{} (3);
     \draw[tree] (3) to node{} (4);
     \draw[tree] (4) to node{} (5);
     \draw[tree] (5) to node{} (6);
     \draw[negedge] (2) to node{} (5);
     \draw[negedge] (1) to node{} (6);

     
     \end{tikzpicture}\\
     (b) $\Phi'$
     \end{minipage}

     
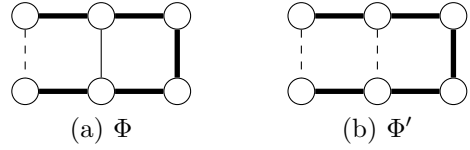
\captionof{figure}{Signed digraphs for Example \ref{example: more than one member of an iso class}.}
     \label{fig: example cospec sw.eq.}
\end{minipage}

\end{table}

The below conclusions then follows immediately.
\begin{corollary}\label{cor: you may choose a spanning tree}
Let $G$ be a graph and let $T\subseteq E(G)$ be a spanning tree on $G$. Further, let $\mathcal{D}$ denote the collection of all signed digraphs on $G$, and let $\mathcal{D}_T\subset\mathcal{D}$ be the collection of such signed digraphs that coincide with $T$ on the relevant  edges. Then  
\[\Phi\in\mathcal{D} ~\iff~ \exists \Phi'\in\mathcal{D}_T \text{ s.t. } \Phi\sim\Phi'.\]
\end{corollary}
The following well-known fact also follows immediately from Proposition \ref{prop: we can always make a tree equal}. 
\begin{corollary}\label{cor: forest}
Let $G$ be a forest and let $\Phi=(G,\varphi)$ be a signed digraph. Then $\Phi\sim G$. 
\end{corollary}
Note that $\mathcal{D}_T$ may contain more than one member from a given switching isomorphism class, as illustrated in the following example. 

\example{example: more than one member of an iso class}{
Consider the non-isomorphic signed digraphs $\Phi$ and $\Phi'$, as illustrated in Figures \ref{fig: example cospec sw.eq.}a and \ref{fig: example cospec sw.eq.}b, respectively. 
It is obvious that a tree, represented by the thick lines, coincides. However, if the vertically oriented  digons in Figure \ref{fig: example cospec sw.eq.}b are multiplied by $-1$ (which clearly is a switching operation), the result is isomorphic to Figure \ref{fig: example cospec sw.eq.}a. 
Thus, $\Phi$ and $\Phi'$ belong to the same switching isomorphism class.
}

As a closing remark to this section, we would like to express interest in the nontrivial question that follows up on the example above, and asks exactly how many members of a given switching equivalence class may coincide in a predetermined spanning tree. 
This matter is not explored further in this work. 

\subsection{Symmetric spectra}
\label{subsec: symmetric}
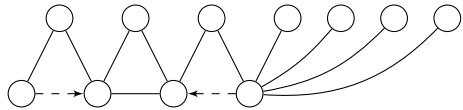
\begin{figure}[h!]
\centering

    \centering
    \begin{tikzpicture}
            \node[vertex] (1) at (0,0) {};
            \node[vertex] (2) at (.5,1) {};
            \node[vertex] (3) at (1,0) {};
            \node[vertex] (4) at (1.5,1) {};
            \node[vertex] (5) at (2,0) {};
            \node[vertex] (6) at (2.5,1) {};
            \node[vertex] (7) at (3,0) {};
            
            \node[vertex] (8) at (3.5,1) {};
            \node[vertex] (9) at (4.2,1) {};
            \node[vertex] (10) at (4.9,1) {};
            \node[vertex] (11) at (5.6,1) {};
            
            \draw[negarc] (1) to node{} (3);
            \draw[edge] (1) to node{} (2);
            \draw[edge] (2) to node{} (3);
            
            \draw[edge] (3) to node{} (4);
            \draw[edge] (4) to node{} (5);
            \draw[edge] (3) to node{} (5);
            
            \draw[negarc] (7) to node{} (5);
            \draw[edge] (6) to node{} (7);
            \draw[edge] (5) to node{} (6);
            
            \draw[edge] (7) to node{} (8);
            \draw[edge,bend right=10] (7) to node{} (9);
            \draw[edge,bend right=18] (7) to node{} (10);
            \draw[edge,bend right=25] (7) to node{} (11);
    \end{tikzpicture}
    \caption{Another signed digraph whose spectrum is symmetric.}
    \label{fig: ugly example}    
\end{figure}

In spectral graph theory, it is commonly asked  which structural characteristics imply symmetry of the corresponding spectrum.
It is well-known that graphs have symmetric spectra if and only if they are bipartite. 
With respect to the (conventional) Hermitian adjacency matrix $H$ \cite{guo2017hermitian}, digraphs have been shown to have symmetric spectra if they are bipartite or (switching equivalent to) an oriented\footnote{A digraph is said to be oriented if it contains no  digons.} digraph, but the reverse implications do not hold. 
In the current context, one may show the following.
\begin{lemma}
Let $G$ be a bipartite graph and let $\Phi=(G,\varphi)$ be a signed digraph. Then the spectrum of $\Phi$ is symmetric around zero.
\end{lemma}
\begin{proof}
Let $\Phi'=(G,\varphi')$, where $\varphi'(uv)=-\varphi(uv)$.
If $G$ is bipartite then every cycle $C$ of $G$ satisfies $\varphi'(C)=(-1)^{|C|}\varphi(C)=\varphi(C)$. 
By Proposition \ref{prop reff}, we thus have $\Phi'\sim \Phi$. 
Finally, since $\E(\Phi')=-\E(\Phi)$ and their spectra coincide, it follows that said spectra are symmetric. 
\end{proof}
However, an oriented signed digraph in general does not have a symmetric spectrum, and no necessary properties were found. 
As a consequence of the existence of sporadic, 'ugly' examples such as the signed digraph in Figure \ref{fig: ugly example},  the authors expect that a tight characterization of all signed digraphs that have symmetric spectra is unlikely to be found. 

%
%

\section{Signed digraphs of low rank}        
\label{sec: low rank}
As we work our way towards spectral characterizations, we first consider signed digraphs of rank $2$ or $3$. 
This restriction severely limits the combinatorial complexity of the corresponding structure. 
Indeed, if almost all rows of its Eisenstein matrix $\E$ need to be linearly dependent, it stands to reason that many vertices have similar relations to one another, as well. 
Without too much effort, one may show the following two results, regarding the ranks of some basic digraphs.
\begin{lemma}
\label{lemma: rank of Pn}
Let $P_n$ be a path of order $n$. 
Then $\Phi=(P_n,\varphi)$ has rank $2\lfloor n/2 \rfloor$ for any $\varphi$.
\end{lemma}
\begin{proof}
Since $\Gamma(\Phi)$ is a tree, $\rank{\Phi}=\rank{\Gamma(\Phi)}=\rank{P_n}$, by Corollary \ref{cor: forest}.
\end{proof}
\begin{lemma}
\label{lemma: rank of Cn}
Let $C_n$ be a cycle of order $n$, with $n$ odd. 
Then $\Phi=(C_n,\varphi)$ has rank $n$ for any $\varphi$.
\end{lemma}
\begin{proof}
$C_n$ has exactly one elementary spanning subgraph, so its characteristic polynomial  $\chi(\lambda)$ has a nonzero coefficient $a_n$, by Theorem \ref{thm: coefficients}, and thus no zero roots. 
\end{proof}
Note that for $n$ even, it also follows that signed digraphs on $C_n$ have rank $n$, except when they have gain $1$ and $n$ is divisible by $4$, or when they have gain $-1$ and $n-2$ is divisible by $4$. 
The case $n=4$ will be relevant later on. 

We note explicitly that since bipartite signed digraphs have symmetric spectra, their ranks are necessarily even. 
Finally, an intuitive observation, concerning the rank of some block matrices, is the following. 
\begin{lemma}\label{lemma: rank 1}
Let $A$ be a Hermitian matrix defined as \[A = \begin{bmatrix}O & A_{12}\\ A_{12}^* & A_{22}.\end{bmatrix},\]
with $A_{12}$ nonzero and where $A_{22}$ has a zero diagonal.
If $\rk A=3$ then $\rk A_{12}=1$ and $2\leq \rk A_{22}\leq 3$. Moreover, if $\rk A = 2$  then $\rk A_{12}=1$ and $A_{22}=O$.
\end{lemma}
Specifically, the lemma above implies that if two vertices are twins in an induced subdigraph of a signed digraph whose rank is $2$ or $3$, then they are also twins in said (larger) signed digraph.

\subsection{Rank 2}
Let us first consider signed digraphs of rank 2. 
As is common in this type of research, the eigenvalue interlacing theorem is used extensively, to forbid particular structures from occurring as induced subgraphs. 
Using the lemmas above, we may characterize the underlying graphs of all signed digraphs with rank 2, based on this idea.
\begin{lemma}
\label{lemma: underlying graph of a rank 2 sdg}
If $\Phi=(G,\varphi)$ is a connected signed digraph with Eisenstein rank $2$, then $G$ is a complete bipartite graph.
\end{lemma}
\begin{proof}
By contradiction. 
By Lemma \ref{lemma: rank of Cn}, $\Phi$ is odd-cycle-free, and thus bipartite. 
Let $P,Q$ denote the coloring classes of $\Phi$. 
Now, suppose to the contrary that $u\in P,~v\in Q$ are a pair of vertices that is nonadjacent in $G$. 
Since $\Phi$ is connected, there is a $u\to v$ path in $G$. 
Let $U\subseteq V(\Phi)$ be the collection of vertices that is traversed on a shortest $u\to v$ path; then $G[U\cup\{u,v\}]\cong P_{2k}$, for some $k>1$.  
However, since $\rk P_{2k} = 2k$, we obtain a contradiction, by interlacing.  
\end{proof}
The next natural question to ask would be which signed digraphs on underlying graph $K_{p,q}$ have Eisenstein rank 2. 
In the below, we show that any signed digraph that satisfies this requirement is switching isomorphic to its underlying graph. 
\begin{proposition}
\label{prop: sdg with rank 2}
If $\Phi=(G,\varphi)$ is a connected signed digraph with Eisenstein rank $2$, then $\Phi$ is switching isomorphic to the complete bipartite graph $G$.
\end{proposition}
\begin{proof}
By Lemma \ref{lemma: underlying graph of a rank 2 sdg}, $G$ is complete bipartite.
Let $P,Q$ denote the coloring classes, as before.  
By Proposition \ref{prop: we can always make a tree equal}, we may without loss of generality choose the  edge gains of a spanning tree of $G$. 
If $u\in P$ and $v\in Q$, then such a spanning tree (say, $T$) may be obtained by taking all  edges incident to at least one of $u$ or $v$. 
If we choose the  edges in $T$ to be positive  digons, the Eisenstein matrix $\E(\Phi)$ contains
\begin{equation}
    \label{eq: proof sdg with rank 2}
    \def\arraystretch{1.2}
    \left[\begin{array}{cc:cc}
     &   &  1 & \mathbf{j}^\top \\
     &    & \mathbf{j} & X \\
      \hdashline
    1 & \mathbf{j}^\top & &  \\
    \mathbf{j} & X^*   & &  \\
    \end{array}\right],
\end{equation}
where the diagonal blocks are square all-zero blocks of appropriate dimensions, $\mathbf{j}$ denotes an all-ones vector and the $X$ blocks are unknown. 
Finally, since all of the induced $4$-cycles must have gain $1$, we have
\[\E(\Phi) = \begin{bmatrix} O_{p\times p} & J_{p \times q} \\ J_{q\times p} & O_{q\times q} \end{bmatrix}.\]
Thus, there is exactly one rank $2$ switching isomorphism class on $K_{p,q}$, and the claim follows. 
\end{proof}
Note that implicitly, all connected rank-$2$ digraphs have underlying graphs that are twin expansions of $K_2$. 
We will see a similar trend if the rank is increased.

\subsection{Rank 3}
Increasing the allowed rank just slightly still allows for a neat characterization of the switching isomorphism classes.
To obtain this characterization, we first obtain an understanding of the twin reduced structure, after which expansions and signatures are included.
\begin{proposition}\label{prop: order 4 rank 3}
Let $\Phi$ be a connected, twin reduced signed digraph of order $4$ and rank $3$. Then 
$\Phi\sim(T_4,\pm),$\footnote{Note that the $(D,\phi)$ notation is used in this instance, as opposed to $(G,\varphi).$} where $T_4$ denotes the order-$4$ transitive tournament.
\end{proposition}
\begin{proof}
Observe that $\Phi$ is not bipartite, as bipartite signed digraphs have even rank. 
Thus, $\Phi$ contains an odd-sized cycle, which implicitly is a triangle. 
Moreover, by connectedness, at least one vertex (say, $s$) in said triangle is also adjacent to the fourth vertex. 
We apply Proposition \ref{prop: we can always make a tree equal} to assume without loss of generality that the  arcs incident to $s$ have gain $1$; that is
\[\E(\Phi) = \begin{bmatrix} 0 & 1 & 1 & 1\\ 1 & 0 & a & \bar{c} \\ 1 & \bar{a} & 0 & b \\ 1 & c & \bar{b} & 0\end{bmatrix},\]
for $a\in\mathbb{T}_6$ and $b,c\in \mathbb{T}_6\cup \{0\}$. (Note that $a,b,$ and $c$ are symmetric, so we may assume w.l.o.g. that $a$ is nonzero.)
Now, note that if $\re a>0$ and $\re b<0$, then by interlacing, $\Phi$ has two positive and two negative eigenvalues, and thus $\rk\E(\Phi) = 4$, and similarly for the pairs $b,c$ and $a,c$. 
Thus, $a,b,c$ must be such that either their real parts are all positive or all negative; we assume positivity for now.  
For $a,b,c$ as above, we then find that 
\begin{equation}\label{eq: determinant rank 3 proof}\det \E(\Phi)=1+|b|+|c| -2\re{ab + ac + bc}.\end{equation}
We proceed to distinguish three cases.
\begin{enumerate}[label=(\roman*)]
    \item $b=c=0$. Then $\det \E(\Phi) = 1$ and $\rank{\Phi}=4$. 
    \item $b\not=0$, $c=0$. Then $\det\E(\Phi) = 2-2\re{ab} = 0 \iff a=\bar{b},$ which implies that $\Phi$ is not twin reduced, contradiction.
    \item $b,c\not=0$.
    Then we need $\det\E(\Phi) = 3-2\re{ab + ac + bc} = 0,$ which holds if and only if one of the following cases is true: (I) $\re{ab} = \re{bc} = \re{ac} = 1/2$ or (II) $\re{ab} = \re{bc} = 1 ~\wedge~ \re{ac}=-1/2.$
    In case (I), we write $ab = \w^{k_a + k_b}$, and use that $\re{ab} = 1/2$ implies  $k_a + k_b$ is odd. 
    Clearly, there are no integers $k_a,k_b,k_c$ such that $(k_a + k_b),(k_b + k_c)$ and $(k_a + k_c)$ are all simultaneously odd, so the desired $a,b,c$ do not exist.
    in case (II), assume $\re a, \re b, \re c > 0$ to conclude that either $a=c=\w$ and $b=\ww$ or $a=c=\ww$ and $b=\w$. Both possibilities are equivalent, and yield $\Phi\sim(T_4,+).$ 
\end{enumerate}
Finally, note that if negativity was assumed, we would have obtained $\Phi\sim(T_4,-)$. 
\end{proof}

Now, we may naturally try to increase the order. A well-known fact is the following.
\begin{lemma}\label{lemma: transitive tournament}
A digraph $D$ is a transitive tournament if and only if every one of its induced subdigraphs is also a transitive tournament. 
\end{lemma}
The above holds analogously when switching is allowed. 
We forego a formal proof, as this fact should be clear by observing that the gains of a basis of the cycle space is known and a quick application of Proposition \ref{prop reff}.

Now, if we extend the above by considering rank-$3$ signed digraphs of order $n\geq 5$, we find that such signed digraphs are not twin reduced. 
\begin{lemma}
\label{lemma: order 5+ rank 3 implies twin}
Let $\Phi$ be a connected signed digraph of order $\geq 5$, rank $3$. Then $\Phi$ is not twin reduced. 
\end{lemma}
\begin{proof}
Let $n\geq 5$ and suppose to the contrary that $\Phi$ is twin reduced, of order $n$ and rank $3$. We distinguish three cases.
\begin{enumerate}[label=(\roman*)]
    \item There is an order-$4$ induced subdigraph $\Phi'$ of $\Phi$ such that $\Gamma(\Phi')\cong K_{1,1,2}.$ This induced subdigraph has rank at most $3$, by interlacing, and is therefore not twin reduced, by Proposition \ref{prop: order 4 rank 3}. Then, using Lemma \ref{lemma: rank 1}, we obtain that the vertices that are twins in $\Phi'$ are also twins in $\Phi$, and thus $\Phi$ is not twin reduced, contradiction.
    \item There is an order-$4$ subdigraph of rank $2$ or $0$. These are respectively complete bipartite or empty, and therefore contain twins; the contradiction follows as in case (i). 
    \item All order-$4$ subdigraphs satisfy Proposition \ref{prop: order 4 rank 3}. Then, by interlacing, they are either all switching isomorphic to $(T_4,+)$ or all to $(T_4,-)$. Hence, using Lemma \ref{lemma: transitive tournament}, we find $\Phi\sim(T_5,\pm)$. However, $\rk{(T_5,\pm)} =5$, and we obtain a contradiction.  
\end{enumerate}
All three possible cases yield a contradiction, so the claim follows.
\end{proof}

The results above show that the number of distinct structures (up to switching isomorphism) whose clique expansions have rank $3$ is very small. Thus, we easily arrive at the following result, that concerns their expanded counterparts.
\begin{proposition}\label{prop: sdg with rank 3}
Let $\Phi$ be a connected signed digraph of rank $3$. Then either $\Phi\sim TE( (K_3,\varphi),\tau)$, for any $\varphi$ and $\tau\in\mathbb{N}^3$, or $\Phi\sim TE((T_4,\pm),\tau)$, $\tau\in\mathbb{N}^4$. 
\end{proposition}
\begin{proof}
Note that any triangle has rank $3$, and observe that any complete tripartite signed digraph has rank $3$ if and only if the  edges between partition groups have equal types, as touched upon in the proof of Lemma \ref{prop: order 4 rank 3}. 
The second part of the claim then follows by application of Lemmas \ref{prop: order 4 rank 3} and \ref{lemma: order 5+ rank 3 implies twin}.
\end{proof}


%
%

\section{Few non-negative eigenvalues}       
\label{sec: extreme ev}
Inspired by the recent work by \citet{oboudi2016characterization}, we inquire into signed digraphs that satisfy a different set of spectral requirements.
Specifically, we now require candidates to have an (almost) minimal number of non-negative or, equivalently, non-positive eigenvalues.  
A neat characterization of signed digraphs that admit to these conditions, could be the foundation upon which to build another set of (possibly infinite) families of signed digraphs, whose spectra determine them, up to switching isomorphism.

Until the end of this section, we will consider the former case, that is, few strictly positive eigenvalues. 
However, note that the negative case is indeed equivalent, and may be obtained by multiplying the signature $\varphi$ by $-1$. 
Additionally, observe that (switching) twins are inherently forbidden, since their existence implies the occurrence of $0$ as an eigenvalue.

\subsection{One non-negative eigenvalue}
\label{sec: one pos ev}
A clear starting point for this investigation is the class of signed digraphs with exactly one non-negative eigenvalue. 
\begin{proposition}\label{prop: second largest eigenv nega implies complete}
The signed digraph $\Phi=(G,\varphi)$ satisfies $\lambda_2<0$ only if $G$ is complete.
\end{proposition}
\begin{proof}
By induction on the order $n$ of $\Phi$. For $n=2,$ the claim is obviously true. 
Now, suppose that the claim holds for any order-$n$ signed digraph $\Phi$, and consider $\Phi'$ of order $n+1$. 
Since the underlying graph of every order-$n$ induced subgraph of $\Phi'$ is complete, by the induction hypothesis and eigenvalue interlacing, it follows that $\Gamma(\Phi')$ is complete. 
\end{proof}

We note explicitly that Proposition \ref{prop: second largest eigenv nega implies complete} is not sufficient.
One may, for example, consider the signed digraph $(T_4,+)$, which has $\lambda_2=0$, while its underlying graph is complete.

Building on Proposition \ref{prop: second largest eigenv nega implies complete}, we will now characterize all signed digraphs whose second largest eigenvalue is negative. 
Indeed, note that we thus far know for sure that such a signed digraph must have a complete underlying graph, but which (and how many) signatures $\varphi$ may be added to yield a signed digraph that satisfies the desired spectral requirements, is as of yet unknown.
In the below, we find that this collection is limited to exactly two switching equivalence classes, for given $n$. 

\begin{definition}\label{def: K*}
Let $K_n^*$ denote the digraph obtained from $K_n$ by orienting exactly one edge. 
\end{definition}
\begin{lemma}\label{lemma: spectrum K*}
$K_n^*$ has a negative second-largest eigenvalue for all $n\geq 3$.  
\end{lemma}
\begin{proof}
The characteristic polynomial of $K_n^*$ is given by 
\[\chi(\lambda)=(\lambda+1)^{n-3}\left(\lambda^3 - (n-3)\lambda^2 -(2n-3)\lambda-1\right),\]
which by Descartes rule of signs has exactly one positive root, and clearly no zero roots.
\end{proof}

\begin{theorem}\label{thm: signed digraph with negative second largest eigenvalue}
Let $\Phi$ be a signed digraph with $\lambda_2<0.$ then either $\Phi\sim K_n$ or $\Phi\sim K_n^*$. 
\end{theorem}

\begin{proof} 
By Proposition \ref{prop: second largest eigenv nega implies complete}, $\Gamma(\Phi)$ is complete. 
Then, the claim may easily be verified for $n\leq 4$, and we proceed by induction. 
Let $n\geq 4$, suppose that the claim is true for signed  digraphs of order $n$ and consider a signed digraph $\Phi$ of order $n+1$. 
Let $k$ denote the number of such induced subdigraphs that are switching isomorphic to $K_n$, and let $V:=V(\Phi')$.

Note that the $n+1$ order-$n$ induced subdigraphs of $\Phi$ must all simultaneously satisfy $\lambda_2<0$, by eigenvalue interlacing. 
This implies that every such order-$n$ subgraph has either only gain-$1$ triangles, or has $n-2$ triangles whose gain is $\omega$, which all intersect on an  edge, and $\binom{n}{3}-n+2$ gain-$1$ triangles.\footnote{Recall that a cycle is said to have gain $\w$ if taking the product of the arc gains corresponding to the arcs hit by traversal of the cycle in at least one direction equals $\w$. I.e., $\varphi(C)=\w$ if either $\varphi(C^\rightarrow)$ or $\varphi(C^\leftarrow)$ equals $\w$.}

Now, we may count the number of pairs $(u,t)$ where $u\in V(\Phi)$ is not part of the gain-$\w$ triangle $t$ in $\Phi$. 
Since for all $u\in V$ we have $\Phi[V\setminus \{u\}]\cong K_n$ or $K_n^*$, where the former holds for $k$ out of the $n+1$ order-$n$ subgraphs, we find that there are $(n+1-k)(n-2)$ such pairs $(u,t)$.
Similarly, since every gain-$\w$ triangle contains all but $n-2$ vertices of $\Phi$, there are $(n-2)\Delta$ such pairs, where $\Delta$ is the total number of gain-$\w$ triangles in $\Phi$. 
Hence, $\Delta=n+1-k$. 

Now, we may distinguish a few cases. 
If $k\geq 4$, then every triangle in $\Phi$ is part of at least one induced subgraph that is switching isomorphic to $K_n$, and thus $T=0$, which implies $k=n+1$ and thus $\Phi\sim K_{n+1}$. 
Similarly, if $k=3$, then all but one triangle in $\Phi$ certainly have gain $1$. 
This implies $\Delta\leq 1$ and thus $n\leq 3$, which is a contradiction. 

If $k=2$, then $\Delta=n-1$. Suppose that $\Phi[V\setminus\{u\}]$ and $\Phi[V\setminus\{v\}]$ are the $\sim K_n$ subdigraphs. 
Then at most the triangles that contain the  edge $(u,v)$ may have gain $\w$; the others all have gain $1$. 
Moreover, since there are exactly $n-1$ such triangles, all of them necessarily have gain $\w$.
Then, using that the collection of triangles in a complete graph form a basis of the cycle space, we may apply Proposition \ref{prop reff} to conclude that $\Phi\sim K_n^*$. 

Next, if $k=1$ then $\Delta=n$. 
Suppose that $\Phi[V\setminus \{u\}]\sim K_n^*$. 
Then, the $n-2$ gain-$\w$ triangles in $\Phi[V\setminus \{u\}]$ intersect on some  edge $(v,z)$. 
Moreover, at least one of $\Phi[V\setminus\{v\}]$ and $\Phi[V\setminus\{z\}]$ is also switching equivalent to $K_n^*$, which thus contains $n-2$ different gain-$\w$ triangles, necessarily containing $u$. 
Hence, $n=\Delta\geq 2(n-2)$, which implies $n\leq 4$ and thus $n=4$. 
This yields the potential counterexample $\Phi\sim K_5^{**}$, where $K_5^{**}$ is obtained from $K_5$ by orienting two of its edges, such that their initial vertex coincides. 
However, a quick computation of its spectrum yields $\lambda_2=0$, and thus $K_5^{**}$ does not satisfy the claim. 

Finally, if $k=0$ then $\Delta=n+1$, which, as above, implies that $n+1\geq 3(n-2)$ and thus $n\leq 3$, which is a contradiction. 
\end{proof}
As was previously mentioned, Theorem \ref{thm: signed digraph with negative second largest eigenvalue} tells us that for a given order $n$, a signed digraph that satisfies $\lambda_2<0$ must belong to exactly one of two (spectrally distinct, recall Lemma \ref{lemma: triangles}) switching isomorphism classes. 
This ties in to a natural spectral characterization result, which is provided in Section \ref{sec: cospectrality}.

To conclude this section, we briefly discuss the natural question how much of the above carries over when instead, signed digraphs with one \textit{positive} eigenvalue are considered; that is, when zero eigenvalues are allowed. 
It turns out that the collection of (twin reduced) signed digraphs with this property contains various ('ugly') members of increasing order, that have little in common with the families of graphs that have been discussed so far.
As such, this is considered out of the scope of this work. 

\subsection{Signed digraphs with $\lambda_2>0>\lambda_3$}
In a recent article, \citet{oboudi2016characterization} characterized all graphs with exactly two non-negative eigenvalues. 
This collection turns out to be an exhaustive list of fairly reasonable length. 
As such, it seems reasonable to ask whether an analogue idea may be applied in the current context.
In this section, we will first find some necessary structural properties, to which any signed digraph that satisfies $\lambda_2>0>\lambda_3$ must admit. 
After that, we will inquire into the signatures of signed digraphs on such graphs. 

\subsubsection{Necessary properties}
The original result by \citet{oboudi2016characterization} follows quite straightforwardly as a forbidden subgraph result that forbids $O_3$ and $C_4$. 
Clearly, $O_3$ should still be forbidden, as its inclusion would imply a non-negative third-largest eigenvalue, by eigenvalue interlacing.
However, since a signed digraph on $C_4$ still meets the requirements if its gain is not $1$, we must substantially deviate from the conclusions in \cite{oboudi2016characterization}. 
As usual, let us first consider the graph structures that may be underlying to signed digraphs that satisfy our needs.

\begin{lemma}\label{lemma: 3k1 free}
Let $G$ be a connected $O_3$-free graph, of order $n\geq 5$.  
Every order-$5$ vertex-induced subgraph contains a $C_5$ or a clique expansion of $P_4$.
\end{lemma}

It should be noted that the collection of graphs that are $O_3$-free contains many graphs with higher edge-density than clique expansions of $C_5$ and $P_4$. However, as should be clear to the reader, given such a graph, one may always remove edges to arrive at a graph that is still $O_3$-free, but which \textit{is} such an expansion. That is to say, a graph is $O_3$-free \textit{because} every relevant subset of the vertices is contained in either one of $C_5$, $P_4$, or a clique. 

Given the above, we may formulate some natural conditions for a signed digraph to satisfy $\lambda_2>0>\lambda_3$. 
These will be particularly useful in Section \ref{sec: cospectrality}, when we are constructing potential cospectral mates of a given signed digraph. 

\begin{proposition}\label{prop: necessary negative third largest ev}
Let $\Phi=(G,\varphi)$ be a connected signed digraph that satisfies $\lambda_2>0>\lambda_3$. Then $G$ is a clique expansion of $P_4$ or $C_5$, possibly supplemented with additional edges up to a complete graph. Additionally, it must satisfy the following:
\begin{enumerate}[label=(\roman*)]
    \item For every $U\subset V(\Phi)$  with $\Gamma(\Phi[U])=K_4$, at least one triangle in $\Phi[U]$ is positive, 
    \item For every $U\subset V(\Phi)$  with $\Gamma(\Phi[U])=C_4$, it holds that $\varphi(\Phi[U])\not=1,$ 
    \item For every $U\subset V(\Phi)$  with $\Gamma(\Phi[U])=C_5$, it holds that $\mathrm{Re}\left(\varphi(\Phi[U])\right)<0,$
\end{enumerate}
\end{proposition}

\begin{proof}
Follows from Lemma \ref{lemma: 3k1 free} and the forbidden subdigraphs switching isomorphic to $(K_4,\varphi_1),$ $(C_4,+)$ and $(C_5,\varphi_2)$, where $\varphi_1$ is such that all triangles in $K_4$ are negative, and $\varphi_2$ is such that the $5$-cycle has positive gain. 
\end{proof}

In case we relax the assumption on connectedness, the following conclusion is an immediate consequence of Theorem \ref{thm: signed digraph with negative second largest eigenvalue}.
\begin{proposition}\label{prop: disjoint cliques l2>0>l3}
    Let $\Phi$ be a signed digraph on $G$, where $G$ is a graph that is obtained as the disjoint union of at least two connected components. If $\Phi$ satisfies $\lambda_2>0>\lambda_3$ then $\Phi = \Phi_1\cup \Phi_2$, where $\Phi_j \sim K_{n_j}$ or $\Phi_j\sim K_{n_j}^*$, $j=1,2$.  
\end{proposition}
The conditions in Proposition \ref{prop: necessary negative third largest ev} are certainly not sufficient; plenty of examples to show this are provided in Figures \ref{fig: problem digraph for semicomplete part}, \ref{fig: example semicomplete neg tri} and \ref{fig: c5 expansions illustrations for proof}, as well as any clique expansions of $C_5$
that exceed Table \ref{tab: C5 lex expac vectors}.

Due to an abundance of possibilities, the full classification of signed digraphs with $\lambda_2>0>\lambda_3$ is not provided here. 
However, we will still zoom in on a few special cases. 
While the complete graph seems like an attractive starting point, the vast number of admissible signatures drove the authors to first consider more palpable families. 
In particular, we will investigate a selection of the clique expansions of $P_4$ and $C_5$, which in a sense have the minimal required number of edges.
In the remainder of this section, we will classify such signed digraphs that satisfy $\lambda_2>0>\lambda_3$; these families will be revisited in Section \ref{sec: cospectrality}, where we provide spectral characterizations.


\subsubsection{Short kite graphs}
An $(a,b)$-kite is said to be obtained from a $K_a$ and a $P_b$ by connecting some vertex in the clique to a pendant vertex of the path. 
Such graphs have recently been shown to be determined by their adjacency spectra, for all choices of $a$ and $b$ \citep{sorgun2016spectral}. 
Moreover, if $b=1$ or $b=2$, then for any $a\geq 2$, the corresponding $(a,b)$-kite graph is $O_3$-free, and may potentially satisfy $\lambda_2>0>\lambda_3$.
In fact, we obtain a rather nice parallel to the results of Section \ref{sec: one pos ev}.
It seems intuitive that the complete part of the kite should have a negative second largest eigenvalue, in order for the corresponding signed digraph to satisfy $\lambda_2>0>\lambda_3$.
An elegant application of eigenvalue interlacing confirms this belief. 

\begin{proposition}\label{prop: complete part of kite}
Let $\Phi = (K_n,\varphi)$, $\mathbf{v}\in \{\mathbb{T}_6\cup 0\}^{n}$ and let 
\[\E := \begin{bmatrix} \E(\Phi) & \mathbf{v} & \mathbf{0} \\ \mathbf{v}^* & 0 & 1 \\ \mathbf{0}^\top & 1 & 0 \end{bmatrix}.\]
Then $\E$ satisfies $\lambda_2>0>\lambda_3$ if and only if $\Phi$ is switching isomorphic to $K_n$ or $K_n^*$.
\end{proposition}
\begin{proof}
Since the eigenvalues $\mu_j$ of \begin{equation}\E':=\begin{bmatrix}\E(\Phi) & \mathbf{0} \\ \mathbf{0}^\top& 0 \end{bmatrix}\nonumber\end{equation} interlace those of of $\E$, necessity of the claim follows by Theorem \ref{thm: signed digraph with negative second largest eigenvalue}.
Indeed, note that $\mu_3\geq 0$ yields a contradiction, by interlacing.

Now suppose that $\Phi$ is switching isomorphic to $K_n$ or $K_n^*$. 
Then, again using that the eigenvalues $\mu_j$ of $\E'$, which by construction  satisfy $\mu_1 > 0 = \mu_2 > \mu_3$, interlace those of $\E$, we obtain  $\lambda_1\geq \mu_1\geq \lambda_2\geq 0 \geq \lambda_3$. Finally, observe that
\[\prod_{j=1}^{n+2}\lambda_j = \det \E =
-\det \begin{bmatrix} \E(\Phi) & \mathbf{0} \\ \mathbf{v}^* & 1\end{bmatrix} =
- \det \E(\Phi) = \mu_1\cdot \prod_{j=3}^{n+1} \mu_j \not=0,\]
and thus $\lambda_2>0>\lambda_3.$
\end{proof}
Since the above is quite independent of the choice of $\mathbf{v}$, the desired result follows easily.
\begin{corollary}\label{cor: kite graphs}
Let $\Phi$ be a signed digraph of order $n\geq 5$ whose underlying graph is a $(n-1,1)$ or $(n-2,2)$-kite graph. Then $\Phi$ satisfies $\lambda_2>0>\lambda_3$ if and only if it contains a subdigraph that is switching isomorphic to $K_{n-2}$ or $K_{n-2}^*$ that is non-adjacent to the pendant vertex.  
\end{corollary}

\subsubsection{Semi-complete signed digraphs}
If not one, but instead both pendants of $P_4$ are expanded, we obtain a convenient structure that contains many induced $(a,2)$-kites. 
Thus, if $\lambda_2>0>\lambda_3$ is required, we get substantial structural information almost for free. 
\begin{definition}
Let $G=CE(P_4,[p~~1~~1~~q]),$ $p,q\in\mathbb{N}$ with $p,q\geq 2$, and let $P\subset V(G)$ be the vertices associated with $p$. 
Let $\tilde{\varphi}$ be the signature that differs from the all-one signature only on $(u,v)$, $u,v\in P$, which has $\tilde{\varphi}(u,v)=\w$. 
Similarly, let $\hat{\varphi}$ be the signature that differs from the all-one signature only on $(s,t)$, $s,t\in Q$, with $\hat{\varphi}(s,t)=\w$.
Then $\tilde{\Phi}:=(G,\tilde{\varphi})$ and $\hat{\Phi}:=(G,\hat{\varphi}).$ 
\end{definition}
\begin{proposition}\label{prop: two cliques and a bridge}
Let $p,q\in\mathbb{N}$, let $G=CE(P_4,[p~~1~~1~~q])$,
Then $\Phi=(G,\varphi)$ satisfies $\lambda_2>0>\lambda_3$ if and only if $\Phi\sim G$, $\Phi\sim\tilde{\Phi}$ or $\Phi\sim\hat{\Phi}$.
\end{proposition}
\begin{proof}
Necessity follows by a straightforward application of Proposition \ref{prop: disjoint cliques l2>0>l3}, while respecting the forbidden subdigraphs in Figure \ref{fig: problem digraph for semicomplete part}.  

Now, suppose that $\Phi\sim G$, $\Phi\sim\tilde{\Phi}$ or $\Phi\sim\hat{\Phi}$ which contains respectively $K_{p+1}\cup K_q$,  $K_{p+1}^*\cup K_q$ or $K_{p}\cup K_{q+1}^*$ as an induced subgraph. Since either satisfies $\mu_2>0>\mu_3$ (by Theorem \ref{thm: signed digraph with negative second largest eigenvalue}), we obtain by interlacing that $\Phi$ has at most three positive eigenvalues. 
Then, using some elementary matrix algebra, we find that both cases satisfy
\[\text{sign}\left(\det \E(\Phi)\right)= (-1)^{p+q+2}.\]
Hence, its number of positive eigenvalues must be even, and thus be equal to two. 
\end{proof}

The $ G$ above is a so-called semi-complete graph, which in general consists of two cliques and an arbitrary number of bridges.
In his investigation of graphs with at most two non-negative eigenvalues, \citet{oboudi2016characterization} finds that graphs, which satisfy $\lambda_2>0>\lambda_3$, are clique expansions the members of a family of clique reduced semi-complete
graphs that are $C_4$-free.

Thus, it would be natural to ask which signed digraphs on semi-complete graphs satisfy $\lambda_2>0>\lambda_3$. 
Conveniently, we find that the clique reduced graphs in Oboudi's family contain large induced $(n-2,2)$-kite graphs, to which we may apply Corollary \ref{cor: kite graphs} and interlacing to conclude that its complete parts should be switching isomorphic\footnote{In case $n$ is even; if $n$ is odd then the situation is slightly more complicated. Here, $K$ and $K^*$ are respectively the complete graph and the signed digraph defined in Definition \ref{def: K*}, of appropriate order.} to $K$ or $K^*$. 
However, this does not yield much useful information regarding the possible signatures of the bridges. 
In fact, it turns out that there are admissible signed digraphs on Oboudi's graphs whose triangle gains do not all share the same sign. 
Moreover, if we generalize from Oboudi's family of graphs, and allow for induced four-cycles, a similar phenomenon occurs. 
These graphs are illustrated in Figure \ref{fig: example semicomplete neg tri}.

These potentially occurring negative triangles open up a vast number of potential signatures to consider.
Thus, a concise, full classification of the signed digraphs with $\lambda_2>0>\lambda_3$, whose underlying graphs are semi-complete graphs may not exist.

\begin{table}[t]
\begin{minipage}{0.48\textwidth}\centering
    \begin{minipage}{.49\textwidth}\centering
     \begin{tikzpicture} [scale = .65]
     \node[vertex] (1) at (0,0) {};
     \node[vertex] (2) at (2,0) {};
     \node[vertex] (3) at (1,1.4142) {};
     \node[vertex] (4) at (1,0.54) {};
     
     \node[vertex] (5) at (0,4) {};
     \node[vertex] (6) at (2,4) {};
     \node[vertex] (7) at (1,4-1.4142) {};
     
     \draw[edge] (1) to node{} (2);
     \draw[edge] (1) to node{} (3);
     \draw[edge] (1) to node{} (4);
     \draw[arc] (2) to node{} (3);
     \draw[edge] (2) to node{} (4);
     \draw[edge] (3) to node{} (4);
     
     \draw[edge] (5) to node{} (6);
     \draw[edge] (5) to node{} (7);
     \draw[edge] (6) to node{} (7);
     
     \draw[edge] (3) to node{} (7);
     
     \end{tikzpicture}\\
     (a)
     \end{minipage}
     \begin{minipage}{.49\textwidth}\centering
     \begin{tikzpicture} [scale = .65]
     \node[vertex] (1) at (0,0) {};
     \node[vertex] (2) at (2,0) {};
     \node[vertex] (3) at (1,1.4142) {};
     
     \node[vertex] (5) at (0,4) {};
     \node[vertex] (6) at (2,4) {};
     \node[vertex] (7) at (1,4-1.4142) {};
     
     \draw[arc] (1) to node{} (2);
     \draw[edge] (1) to node{} (3);
     \draw[edge] (2) to node{} (3);
     
     \draw[arc] (5) to node{} (6);
     \draw[edge] (5) to node{} (7);
     \draw[edge] (6) to node{} (7);
     
     \draw[edge] (3) to node{} (7);
     
     \end{tikzpicture}\\
     (b)
     \end{minipage}

     \captionof{figure}{Two signed digraphs with $\lambda_3=0.$}
     \label{fig: problem digraph for semicomplete part}
\end{minipage}~~~~
\begin{minipage}{0.48\textwidth}\centering
    \begin{minipage}{.49\textwidth}\centering
     \begin{tikzpicture} [scale = .65]
     \node[vertex] (1) at (0,0) {};
     \node[vertex] (2) at (2,0) {};
     \node[vertex] (3) at (1,1.4142) {};
     
     \node[vertex] (5) at (0,4) {};
     \node[vertex] (6) at (2,4) {};
     \node[vertex] (7) at (1,4-1.4142) {};
     
     \draw[edge] (1) to node{} (2);
     \draw[edge] (1) to node{} (3);
     \draw[edge] (2) to node{} (3);
     
     \draw[edge] (5) to node{} (6);
     \draw[edge] (5) to node{} (7);
     \draw[edge] (6) to node{} (7);
     \draw[edge,bend right=15] (3) to node{} (6);
     \draw[arc] (3) to node{} (7);
     \draw[negedge, bend right=15] (2) to node{} (6);
     
     \end{tikzpicture}\\
     (a) 
     \end{minipage}
     \begin{minipage}{.49\textwidth}\centering
          \begin{tikzpicture} [scale = .65]
     \node[vertex] (1) at (0,0) {};
     \node[vertex] (2) at (2,0) {};
     \node[vertex] (3) at (1,1.4142) {};
     
     \node[vertex] (5) at (0,4) {};
     \node[vertex] (6) at (2,4) {};
     \node[vertex] (7) at (1,4-1.4142) {};
     
     \draw[edge] (1) to node{} (2);
     \draw[negarc] (3) to node{} (1);
     \draw[edge] (2) to node{} (3);
     
     \draw[edge] (5) to node{} (6);
     \draw[arc] (7) to node{} (5);
     \draw[edge] (6) to node{} (7);
     \draw[arc] (5) to node{} (1);
     \draw[arc] (3) to node{} (7);
     \draw[edge] (2) to node{} (6);
     
     \end{tikzpicture}\\
     (b) 
     \end{minipage}

     \captionof{figure}{Two semi-complete examples with $\lambda_2>0>\lambda_3$, which contain negative triangles.}
     \label{fig: example semicomplete neg tri}
\end{minipage}
\end{table}

\subsubsection{Clique expansions of $C_5$}\label{sec: cycles}
By applying what we have learned about kite graphs, we may draw some interesting conclusions with regard to the expansions of $C_5$.
Indeed, it is not hard to see that any clique expansion of $C_5$ contains many proper induced subgraphs that are simply $(n-2,2)$-kite graphs. 
Thus, under the usual assumptions (recall Prop. \ref{prop: necessary negative third largest ev}), we may substantially limit the potential signatures.

\begin{figure}[b]
     \begin{center}
     \begin{subfigure}[b]{0.20\textwidth}\centering
     \begin{tikzpicture} [scale=1.3]
     
    \node[group] (g1) at (1,0){$K_{t1}$};
    \node[group] (g2) at (0,0){$K_{t2}$};
    \node[group] (g3) at (1,1){$K_{t3}$};
    \node[group] (g4) at (.5,1.8){$K_{t4}$};
    \node[group] (g5) at (0,1){$K_{t5}$};
    
    \draw[edge] (g1) to node{} (g2);
    \draw[edge] (g2) to node{} (g3);
    \draw[edge] (g3) to node{} (g4);
    \draw[edge] (g4) to node{} (g5);
    \draw[negedge] (g5) to node{} (g1);
     
     \end{tikzpicture}
     \caption{}
     \label{fig: 1}
     
     \end{subfigure}
     \begin{subfigure}[b]{0.20\textwidth}\centering
     \begin{tikzpicture} [scale=1.3]
     
    \node[group] (g1) at (1,0){$K_{t1}$};
    \node[group] (g2) at (0,0){$K_{t2}$};
    \node[group] (g3) at (1,1){$K_{t3}$};
    \node[group] (g4) at (.5,1.8){$K_{t4}$};
    \node[group] (g5) at (0,1){$K^*_{t5}$};
    
    \draw[edge] (g1) to node{} (g2);
    \draw[edge] (g2) to node{} (g3);
    \draw[edge] (g3) to node{} (g4);
    \draw[edge] (g4) to node{} (g5);
    \draw[negedge] (g5) to node{} (g1);
     
     \end{tikzpicture}
     \caption{
     }
     \label{fig: 2}
     
     \end{subfigure}
     \begin{subfigure}[b]{0.20\textwidth}\centering
     \begin{tikzpicture} [scale=1.3]
     
    \node[group] (g1) at (1,0){$K_{t1}$};
    \node[group] (g2) at (0,0){$K_{t2}$};
    \node[group] (g3) at (1,1){$K_{t3}$};
    \node[group] (g4) at (.5,1.8){$K_{t4}$};
    \node[group] (g5) at (0,1){$K_{t5}$};
    
    \draw[edge] (g1) to node{} (g2);
    \draw[edge] (g2) to node{} (g3);
    \draw[edge] (g3) to node{} (g4);
    \draw[edge] (g4) to node{} (g5);
    \draw[negarc] (g5) to node{} (g1);
     
     \end{tikzpicture}
     \caption{}
     \label{fig: 3}
     
     \end{subfigure}
     \begin{subfigure}[b]{0.37\textwidth}\centering
     \begin{tikzpicture} [scale=1.2]
     
    \node[shape=circle,draw,dotted,minimum size = 4em] (g1) at (0,0){};
    \node[group] (g2) at (1.5,0){$K_{t2}$};
    \node[group] (g3) at (2.1,1.5){$K_{t3}$};
    \node[group] (g4) at (.75,1.5){$K_{t4}$};
    \node[shape=circle,draw,dotted,minimum size = 4em] (g5) at (-.6,1.5){};
    \draw (-0.8,0) node {$K_{t1}$};
    \draw (-1.5,1.5) node {$K_{t5}$};
    
    \node[vertex] (g1n1) at (-.25,-.15) {};
    \node[vertex] (g1n2) at (.25,.15) {};    
    \node[vertex] (g5n1) at (-.85,1.35) {};
    \node[vertex] (g5n2) at (-.35,1.65) {};
    
    \draw[edge] (g1) to node{} (g2);
    \draw[edge] (g2) to node{} (g3);
    \draw[edge] (g3) to node{} (g4);
    \draw[edge] (g4) to node{} (g5);
    
    \draw[negarc] (g5n1) to (g1n1);
    \draw[negedge] (g5n1) to (g1n2);
    \draw[negedge] (g5n2) to (g1n1);
    \draw[negedge] (g5n2) to (g1n2);
    
    \draw[edge] (g1n1) to (g1n2);
    \draw[edge] (g5n1) to (g5n2);
     
     \end{tikzpicture}
     \caption{
     }
     \label{fig: 4}
     
     \end{subfigure}
     
     \end{center}
     \caption{\label{fig: stomme signature types want hoe moet ik dit anders opschrijven} Illustrations of signature types A, B, C, D. 
     }
     \end{figure}
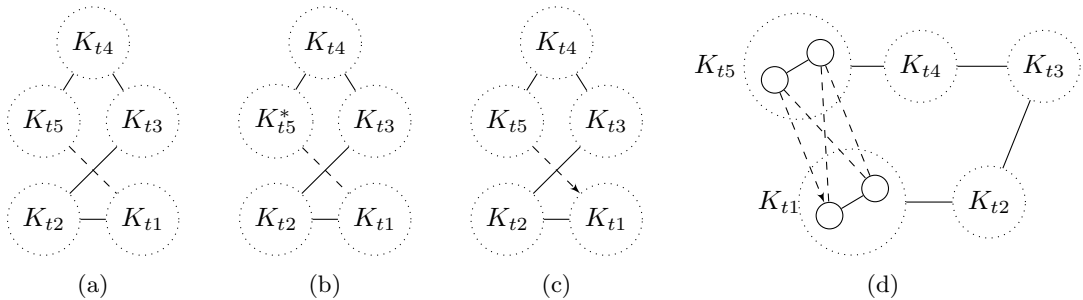

In order to structure the following discussion, it is convenient to define some distinct types of signatures, for clique expansions of $C_5$. 
Given some $ G=CE(C_5,\tau)$ and $\Phi=( G,\varphi)$, the four distinct signature types $\varphi$ are displayed in Figure \ref{fig: stomme signature types want hoe moet ik dit anders opschrijven}. 
Informally, all induced cliques in types $A$ and  $C$ are switching isomorphic to complete graphs, while all of their induced $5$-cycles have gains $-1$ and $-\w$, respectively. 
Oppositely, types $B$ and $D$ do contain gain-$\w$ triangles. 
If the induced cliques associated with expansion parameter $\tau_j$ are denoted $C_j$, then type $D$ is such that exactly one pair $(i,j)$ is such that $C_i\cup C_j$ induces a $K^*$, while every $C_j$ induces $K$ and the induced $5$-cycles have gain $1$ or $\omega$.  
Similarly, expansions of type $B$ are such that exactly one $C_j$ induces a $K^*$, the remaining four induce $K$, and all induced $5$-cycles have gain $1$. 
%
We find the following. 

\begin{proposition}\label{prop: expac of C5 with l2>0>l3}
Let $ G$ be a clique expansion of $C_5$, and let $\Phi=( G,\varphi)$ be a signed digraph that satisfies $\lambda_2>0>\lambda_3$. Then $\Phi\sim\Phi'=( G,\varphi')$ where $\varphi'$  is type A, B, C or D.
\end{proposition}
\begin{proof}
As usual, we may assume that a spanning tree $Y\subseteq E( G)$ of the  edges are positive digons in $\Phi$. 
Specifically, if we denote the cliques corresponding to expansion coefficient $t_j$ (see Figure \ref{fig: stomme signature types want hoe moet ik dit anders opschrijven}) by $ G_j$, a convenient choice of spanning tree is obtained by fixing five nodes $u_j\in V( G_j)$, $j=1,\ldots,5$, and choosing the spanning tree 
\[Y = \bigcup_{j=1}^5\left\{(u_j,v) ~|~ v\in V( G_j)\setminus u_j \right\} \cup \bigcup_{j=1}^4\{(u_j,u_{j+1})\}.\]
Since the subgraph $\Phi[V( G_j)\cup V( G_k)]$ induced by two adjacent cliques, indexed by $j,k$, is again a clique, the subgraph $\Phi[V( G_j)\cup V( G_k)]$ necessarily satisfies $\mu_2<0$, since the eigenvalues of $\Phi[V( G_j)\cup V( G_k)]\cup O_1$ interlace those of $\Phi$. 
For such subgraphs with $k=j+1$, $j\in[4]$,  
we may use that a spanning tree of the induced clique consists of positive  digons, to find that exactly one of two cases must be true: either all  edges in the induced clique are positive  digons, or the induced clique contains exactly one positive  arc and the remainder is made up of positive  digons. 
Finally, $\Phi[V( G_1)\cup V( G_5)]$ must also be switching equivalent to either $K_m$ or $K_m^*$, for appropriate $m$. However, since every induced $C_5$ necessarily has negative gain, it follows straightforwardly that the  edges between $ G_1$ and $ G_5$ must either be all negative  digons, or a single negative  arc supplemented with negative  digons.
Note that either option may indeed be obtained from $K_m$ or $K_m^*$ with a simple diagonal switch that hits the  edges between $ G_1$ and $ G_5.$

By the above, no two adjacent $ G_j$ may both contain an  arc.
Thus, natural question would be whether or not two non-adjacent cliques may both contain an  arc. 
Now, since the smallest admissible signed digraph that might satisfy this property, structured as Figure \ref{fig: two K3*},
has a zero third largest eigenvalue, we may conclusively answer this question with 'no.'
In a similar vein, since the structure in Figure \ref{fig: neg arc and K3*} also has $\lambda_3=0$, none of the cliques may contain an  arc if one of the induced five-cycles has gain $-\omega$. 
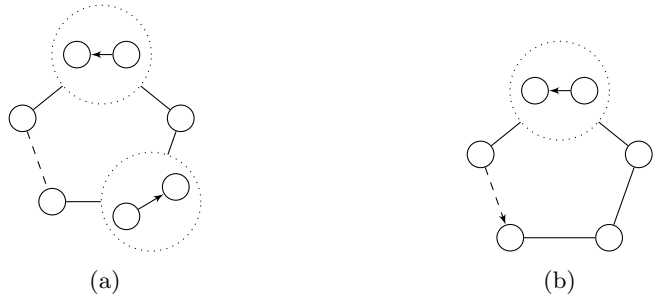
\begin{figure}[t]
    \centering
    \begin{subfigure}[b]{0.40\textwidth}\centering
    \begin{tikzpicture}[scale=1.3]
             \node[vertex] (1) at (0,0) {};
             \node[group] (2) at (1,0) {~~~~~~~~~~};
             \node[vertex] (3) at (1.3,.85) {};
             \node[group] (4) at (.5,1.5) {~~~~~~~~~~};
             \node[vertex] (5) at (-.3,.85) {};
             \node[vertex] (2a) at ( .75,-.15 ) {};
             \node[vertex] (2b) at ( 1.25,.15 ) {};
             \node[vertex] (4a) at ( .75, 1.5) {};
             \node[vertex] (4b) at (.25 ,1.5 ) {};
             
             \draw[edge] (1) to node{} (2);
             \draw[edge] (2) to node{} (3);
             \draw[edge] (3) to node{} (4);
             \draw[edge] (4) to node{} (5);
             \draw[negedge] (5) to node{} (1);
             
             \draw[arc] (2a) to node{} (2b);
             \draw[arc] (4a) to node{} (4b);
    \end{tikzpicture}
    \caption{}
    \label{fig: two K3*}
     
     \end{subfigure}
     \begin{subfigure}[b]{0.40\textwidth}\centering
    
    \begin{tikzpicture}[scale=1.3]
             \node[vertex] (1) at (0,0) {};
             \node[vertex] (2) at (1,0) {};
             \node[vertex] (3) at (1.3,.85) {};
             \node[group] (4) at (.5,1.5) {~~~~~~~~~~};
             \node[vertex] (5) at (-.3,.85) {};
             \node[vertex] (4a) at ( .75, 1.5) {};
             \node[vertex] (4b) at (.25 ,1.5 ) {};
             
             \draw[edge] (1) to node{} (2);
             \draw[edge] (2) to node{} (3);
             \draw[edge] (3) to node{} (4);
             \draw[edge] (4) to node{} (5);
             \draw[negarc] (5) to node{} (1);
             
             \draw[arc] (4a) to node{} (4b);
    \end{tikzpicture}
    \caption{}
    \label{fig: neg arc and K3*}
     \end{subfigure}
    \caption{Two clique expansions of $C_5$ with $\lambda_3=0$}
    \label{fig: c5 expansions illustrations for proof}
\end{figure}
Combining all of the above, we obtain that if $ G=CE(C_5,\tau)$ and $\Phi$ satisfies $\lambda_2>0>\lambda_3$, then $\varphi$ must be switching equivalent to a type A, B, C, or D  signature. 
\end{proof}
However, as was briefly mentioned before, not any clique expansion of $C_5$ may be underlying to a signed digraph that fits our requirements. 
Using our knowledge on the admissible signatures, we learn the following by a computer search. 
In the below, we write $\mathcal{T}_j = \left\{  G\cong  G'[U] ~|~ U\subseteq V( G'),~  G'=CE(C_5,\tau^{j})\right\};$ that is, $\mathcal{T}_j$ is the collection of all graphs that are obtained from $C_5$ by clique expansion with expansion vector at most $\tau^{j}$. 
\begin{proposition}
\label{prop: expansions of C5: sufficiency part}
Let $\Phi = ( G,\varphi)$ be a signed digraph that satisfies $\lambda_2>0>\lambda_3$. Then, up to switching equivalence, the following holds:
\begin{enumerate}[label=(\roman*)]
\item if $\varphi$ is type A, then $ G\in \bigcup_{j=1}^{13}\mathcal{T}_j$,
\item if $\varphi$ is type B or D, then $ G\in\mathcal{T}_{12}\cup \mathcal{T}_{13},$ and
\item if $\varphi$ is type C, then $ G\in\mathcal{T}_{14},$ 
\end{enumerate}
where the $\tau_j$ are displayed in Table \ref{tab: C5 lex expac vectors}.
\end{proposition}
\begin{table}[t]
\centering
\begin{tabular}{@{}lcccccccccccccc@{}}
\toprule
 $\varphi$ type:     & \multicolumn{11}{c}{A}                                                                                                  &  \multicolumn{2}{c}{A,B,D}                   &  C      \\
      & $\tau^{1}$ & $\tau^{2}$ & $\tau^{3}$ & $\tau^{4}$ & $\tau^{5}$ & $\tau^{6}$ & $\tau^{7}$ & $\tau^{8}$ & $\tau^{9}$ & $\tau^{10}$ & $\tau^{11}$ & $\tau^{12}$ & $\tau^{13}$ &$\tau^{14}$   \\ \cmidrule(r){2-12} \cmidrule(lr){13-14} \cmidrule(lr){15-15} 
$t_1$ & 3        & 3        & 3        & 3        & 4        & 5        & 5        & 5        & 3        & $t_1$       & $t_1$       & $t_1$  & $t_1$           & $t_1$                \\
$t_2$ & 3        & 3        & 4        & 2        & 2        & 3        & 2        & 1        & 1        & 1           & $t_2$       & $t_2$ & 1                & 1                   \\
$t_3$ & 3        & 2        & 2        & 4        & 2        & 1        & 2        & 3        & 5        & 2           & 2           & 1     & 1              & 1                 \\
$t_4$ & 2        & 2        & 2        & 2        & 2        & 2        & 2        & 2        & 2        & 2           & 1           & 1     & $t_4$             & 1                \\
$t_5$ & 1        & 2        & 1        & 1        & 2        & 1        & 1        & 1        & 1        & 1           & 2           & $t_5$     & $t_5$          & 1                       \\ \bottomrule
\end{tabular}
\caption{Maximum clique expansions sizes of $C_5$, such that it admits a signed graph that satisfies $\lambda_2>0>\lambda_3$. Free variables in $\tau^{10}$-$\tau^{14}$ may be arbitrarily large. Note that for types $A$ and $D$, the collections $\mathcal{T}_{12}$ and $\mathcal{T}_{13}$ are equal.
However, since they do not necessarily coincide for type $B$, the distinction is kept.}
\label{tab: C5 lex expac vectors}
\end{table}     
From Table \ref{tab: C5 lex expac vectors}, we may observe that there are still arbitrarily large expansions of $C_5$ that satisfy our needs, in addition to some subtly differently structured smaller ones. 

It should be noted here, that the signed digraphs obtained by taking a $ G\in\mathcal{T}_{j}$ and an admissible $\varphi$ from the parameters and structures described above are, as has been the habit throughout, in some sense leading members of a switching equivalence class. 
For example, if one starts with a $C_5$ whose gain is $-\omega$, any single vertex may be clique expanded to arbitrary size, without compromising the spectral requirement. However, since every signed digraph obtained in such a manner is switching equivalent to one obtained by expansion of vertex "1," as in Figure \ref{fig: stomme signature types want hoe moet ik dit anders opschrijven}, these are not explicitly listed.


%
%

\section{Cospectrality and determination}    
\label{sec: cospectrality}
Of particular interest to the authors are uniquely occurring spectra of graphs. 
That is, spectra that uniquely determine a graph, up to isomorphism.
This notion has received considerable attention for several decades \cite{vandam2003, vandam2009}. 
When an analogous line of research was launched for the Hermitian adjacency matrix $H$, \citet{mohar2016Hermitian} subtly shifted the definition of "determined by the spectrum" such that 'DS' digraphs were now allowed to have non-isomorphic cospectral mates, as long as they were all \textit{switching} isomorphic. 

In a previous work \cite{wissing2019negative}, the authors considered the traditional notion, applied in the Hermitian adjacency matrix paradigm. 
It was determined that digraphs whose $H$-spectra occur uniquely, up to isomorphism, are extremely rare; though some infinite families do exist. 
However, we find that with respect to $\E$, any non-empty signed digraph has a (in fact, many) non-isomorphic, switching equivalent partner. 
This is formally shown below, with an intuitive counting argument. 
Thereafter, we make use of the classifications from Sections \ref{sec: low rank} and \ref{sec: extreme ev} to prove that several of the discussed families have spectra that occur only for their respective switching equivalence classes; i.e., these families are determined by their spectra in the broader sense \cite{mohar2016Hermitian}.

\subsection{Existence of a switching equivalent partner}
\label{subsec: strong determination}

We will formally show that any non-empty signed digraph has at least one switching equivalent, non-isomorphic partner. 
In the upcoming proofs, let $\Omega_k(\E)$ denote the number of entries of $\E$ that are equal to $\omega^k$, for $k\in\mathbb{Z}_6$.
We make the following observation.

\begin{lemma}\label{lemma: cutset isomorphic partner}
Let $\Phi=( G,\varphi)$ be a signed digraph.
Let $U \subset V$ and $W = V\setminus U$ be a cut, and partition $\E(\Phi)$ such that
\[\E = \begin{bmatrix} \E_{U,U} & \E_{U,W} \\ \E_{W,U}& \E_{W, W}\end{bmatrix}.\] 
If there are $k,l\in\mathbb{Z}_6$ such that $\Omega_k(\E_{U,W}) \not = \Omega_l(\E_{U,W})$, then  there is a $\Phi'\not\cong\Phi$ such that $\Phi'\sim\Phi$. 
\end{lemma}

\begin{proof}
Suppose that there are $k,l\in\mathbb{Z}_6$ such that $\Omega_k(\E_{U,W}) \not = \Omega_l(\E_{U,W})$ and assume to the contrary that $\Phi'\cong\Phi$ for all $\Phi'\sim\Phi$. Consider the switching matrix 
\[S_k = \begin{bmatrix}\omega^k I_u & O \\ O & I_{n-u}\end{bmatrix} ~~\text{ with } k\in\mathbb{Z}_6,\]
and set $\E_k=S_k\E S_k^{-1},$ which serves as the Eisenstein matrix of the switching isomorphic signed digraph $\Phi_k$.   
Since isomorphic signed digraphs necessarily contain an equal number of positive  digons, $\Omega_0(\E_k)=\Omega_0(\E)$.
Note that $\E_{U,U}=(\E_k)_{U,U}$ and $\E_{W,W}=(\E_k)_{W,W}$, and that
\[\Omega_{p}(\E_{U,W}) =\Omega_{p+k}((\E_k)_{U,W})~~\text{ for all }p,k\in\mathbb{Z}_6.\]
Then, since $\E,\E_k$ are Hermitian, it follows that $\Omega_0(\E_{U,W})=\Omega_0(\E_k)_{U,W}=\Omega_{-k}(\E_{U,W})$ for every $k\in\mathbb{Z}_6$ and we obtain a contradiction.  
\end{proof}

Now, we may simply consider the number of  edges that is hit by a given cut, to determine that the required cut $U$ and $k,l$ certainly exist in a given non-empty signed digraph. 
\begin{proposition}
Let $\Phi$ be non-empty. Then there exists a $\Phi'\not\cong\Phi$ such that $\Phi\sim\Phi'$.
\end{proposition}
\begin{proof}
Let $\Phi$ be a non-empty signed digraph, and suppose that $\Phi$ is strongly determined by its spectrum.  
Then, for any cut $U\subset V(D)$ and any switching over the  edges between $U$ and $V(\Phi)\setminus U$, the digraph obtained by the corresponding switching is isomorphic to $\Phi$. 
By Lemma \ref{lemma: cutset isomorphic partner}, this implies that any cut of $\Phi$ hits equally many  edges of every type.
This, in turn, implies that any cut in $\Gamma(\Phi)$ must hit a number of  edges that is divisible by $6$. 

Now, suppose that $u,v$ are two vertices that are neighbors in $\Gamma(\Phi)$. By the above, the degrees $d_u$ and $d_v$ of $u$ and $v$, respectively, must satisfy $d_u \equiv 0~(\text{mod}~6)$ and $d_v \equiv 0~(\text{mod}~6).$ But then the cut set $\{u,v\}$ hits exactly $d_u + d_v - 2 \equiv 4~(\text{mod}~6)$  edges, and we have a contradiction.
\end{proof}
To conclude, we remark that the line above holds for any finite unit gain group except $\mathbb{T}_2$. 
However, if the the set of allowed complex unit gains is not closed under multiplication, one may construct examples such that each cut that violates the premise above does not allow for gain switching, and therefore cannot produce a counterexample. 
This phenomenon occurs, for example, with the Hermitian adjacency matrix for directed graphs, and is exploited in \cite{wissing2019negative}.


\subsection{Spectral determination}
By the conclusion in the previous section, it is natural to adopt the before mentioned definition of spectral determination due to Mohar. 
Formally, we have the following. 

\begin{definition}{\cite{mohar2016Hermitian}}
Let $\Phi$ be a signed digraph and let $Co(\Phi)
$ be the collection of signed digraphs whose spectra coincide with the spectrum of $\Phi$. 
$\Phi$ is said to be \textit{determined by its $\E$-spectrum} (DES) if $\Phi\sim\Phi'$ for all $\Phi'\in Co(\Phi)$.
\end{definition}

In the remainder of this work, we draw from the classifications in Sections \ref{sec: low rank} and \ref{sec: extreme ev} and verify whether or not some of these families, whose spectral behaviour is, in a sense, extreme, have non-equivalent cospectral mates. 

\subsubsection{Low rank}
\label{subsec: weak determination}
Let us first consider the families of signed digraphs with low rank. 
Since the collection of graphs that might be underlying to such signed digraphs was neatly characterized, we may straightforwardly show the following results. 
\begin{proposition}\label{prop: rank 2 DES under connectedness}
Let $\Phi$ be a connected signed digraph with $\rank{\Phi}=2$. If $D$ is connected and cospectral to $\Phi$, then $D\sim \Phi$. 
\end{proposition}
\begin{proof}
Recall from Proposition \ref{prop: sdg with rank 2} that $\rank{D}=\rank{\Phi}=2$ implies $\Phi\sim K_{f,g}$ and $D\sim K_{p,q}$, for $f,g,p,q\in \mathbb{N}$.
Now, we may simply solve 
\begin{align*}\begin{cases} |V(\Phi)| = |V(D)| \\ |E(\Phi)| = |E(D)|
\end{cases} \iff 
    \begin{cases} p + q = f + g\\ pq = fg\end{cases} & 
    \iff (p,q) = (f,g) ~\vee (p,q) = (g,f).
\end{align*}
To obtain $\Phi\sim K_{f,g}\cong K_{p,q}\sim D$, which completes the proof.
\end{proof}
An important note to place here is that the assumption on connectedness is almost always required. Indeed, note that for instance $K_{1,4}$ and $K_{2,2}\cup K_1$ (known as the saltire pair) admit to the requirements, but are cospectral. 
The reason is quite simple: if connectedness is relaxed, then one may (using the notation from the proof above) simply find numbers $p,q$ such that $pq=fg$, add $r$ isolated vertices to satisfy $p+q+r=f+g$. 
The following small generalizations follow straightforwardly from this insight.
\begin{corollary}\label{cor: rank 2 DES if prime}
Let $\Phi\sim K_{p,q}$ with $p,q$ prime. Then $\Phi$ is DES. 
\end{corollary}
\begin{corollary}
Let $\Phi\sim K_{n,n}$, for $n\in\mathbb{N}.$ Then $\Phi$ is DES.
\end{corollary}
\begin{proof}
Follows since $\Phi$ attains the minimum number of vertices ($2n$) necessary for a rank-$2$ signed digraph with $|E|=n^2$. 
\end{proof}
In Section \ref{sec: low rank}, we have shown that $\Phi$ has rank $3$ if and only if its twin reduction is equivalent to either a triangle\footnote{Recall that there are exactly four classes of triangles, which contain $K_3$, $K_3^*$, $-K_3$ and $-K_3^*$, respectively.} or an order-$4$ transitive tournament.
However, contrary to the above, their respective expansions may sometimes be cospectral, as was observed in \cite{li2021hermitian}. 
Cospectrality occurs for each pairing of the positive reduced graphs $K_3, K_3^*, T_4$,  and equivalently for their negative counterparts. 
The smallest examples to this fact are:
\begin{itemize}
    \item $TE(K_3,[1~8~15])$ is cospectral to $TE(K^*_3,[3~5~16])$,
    \item $TE(K_3^*,[3~4~7])$ is cospectral to $TE(T_4,[1~1~6~6])$,
    \item $TE(K_3,[3~20~25])$ is cospectral to $TE(T_4,[3~5~10~30])$.
\end{itemize}
Moreover, each of the three reduced signed digraphs has an arbitrarily large  (nontrivial) twin expansion cospectral to a signed digraph with a different underlying graph.
The propositions below are easily checked by simply counting the number of vertices, edges, and triangles\footnote{Since their gains are a factor, the triangles in expansions of $K_3$ should all be counted twice.}. 
\begin{proposition}
\label{lemma: large cospectral rank 3 k3* k3}
Let $\Phi$ and $\Phi'$ be defined as $\Phi = TE(K_3^*,[2i+1~~i(3i+2)~~2(3i+1)(i+1)])$ and $\Phi'=TE(K_3,[i~~(3i+1)(i+1)~~ 2(3i+1)(i+1)-1])$ for $i\in\mathbb{N}.$ Then $\Phi$ is cospectral to $\Phi'$. 
\end{proposition}
\begin{proposition}
\label{lemma: large cospectral rank 3 k3* t4}
Let $\Phi$ and $\Phi'$ be defined as $\Phi = TE(K_3^*,[i(i+1)/2~~i(i+1)/2+1~~i(i+1)+1])$ and $\Phi'=TE(T_4,[1~~i(i-1)/2~~(i+1)(i+2)/2~~i(i+1)])$ for $i\in\mathbb{N}.$ Then $\Phi$ is cospectral to $\Phi'$. 
\end{proposition}

\subsubsection{A single non-negative eigenvalue}
Since the collection of signed digraphs that satisfy $\lambda_1>0>\lambda_2$ on a given number of vertices is characterized by just two switching equivalence classes, we may use the structural information obtained in Section \ref{sec: one pos ev} to draw some quick conclusions with regard to their cospectrality.
\begin{proposition}
Let $\Phi$ be either $K_n$ or $K_n^*$. Then $\Phi$ is DES.
\end{proposition}
\begin{proof}
Let $D$ be cospectral to $\Phi$. Then $D$ has $\lambda_1>0>\lambda_2$, and thus by Theorem \ref{thm: signed digraph with negative second largest eigenvalue}, either $D\sim K_n$ or $D\sim K_n^*$.
Suppose w.l.o.g. that $\Phi=K_n$ and $D\sim K_n^*$. Then, by Lemma \ref{lemma: triangles} or \ref{lemma: spectrum K*}, $D$ is not cospectral to $\Phi$, which is a contradiction. Thus, $D\sim \Phi$. 
\end{proof}
Naturally, the same argument also holds when all  edge gains are multiplied by $-1$. 
\begin{corollary}
Let $\Phi$ be either $(K_n,-)$ or $(K_n^*,-)$. Then $\Phi$ is DES.
\end{corollary}

\subsubsection{Smallest eigenvalue $-1$}
A straightforward example to show that DES signed digraphs do not necessarily consist of a single connected component, possibly appended with a collection of disjoint vertices, carries over from graph theory. 
\begin{lemma}\label{lemma: smallest ev -1}
Let $\Phi$ be an order-$n$ connected signed digraph with $\lambda_n=-1$. Then $\Phi\sim (K_n,+)$.
\end{lemma}
\begin{proof}
Since $\Phi$ is connected, $\Gamma(\Phi)$ contains a shortest path $P_{u,v}$ for every $u,v\in V(\Phi)$. 
The result follows by two applications of interlacing: first to see that every such $P_{u,v}$ is length at most two, and thus $\Gamma(\Phi)=K_n$, and then to obtain that every induced triangle has gain $1$.
\end{proof}
\begin{proposition}\label{prop: smallest ev -1 disjoint sequence}
Let $m\in\mathbb{N}$ and  $\Phi=\bigcup_{j=1}^m (K_{n_j},+)$,  $n_j\in\mathbb{N}, ~j\in[m].$ Then $\Phi$ is DES.
\end{proposition}
\begin{proof}
Note that $\Phi$ has smallest eigenvalue $-1$, and let $D$ be cospectral to $\Phi$. 
Then by Lemma \ref{lemma: smallest ev -1}, every connected component of $D$ is switching isomorphic to a complete graph of appropriate order. 
The result follows since every positive eigenvalue characterizes such a connected component, and every zero eigenvalue corresponds to an isolated vertex.  
\end{proof}
\begin{corollary}
Let $m\in\mathbb{N}$ and  $\Phi=\bigcup_{j=1}^m (K_{n_j},-)$,  $n_j\in\mathbb{N}, ~j\in[m].$ Then $\Phi$ is DES.
\end{corollary}
\subsubsection{Weakly determined with $\lambda_2>0>\lambda_3$}
\label{subsubsec: wds with l3<0}
In the final section of this work, we will draw from the families of graphs characterized in Section \ref{sec: extreme ev}, and use much of their inherent structure to obtain several more families of signed digraphs that are weakly determined by their $\E$-spectra. 

In the below, we will use the same line of proof in two distinct situations, mainly separated by the numbers of edges and triangles that must (at least) occur in the corresponding underlying graphs. 
The considered graphs consistently contain relatively large cliques, which translate to a large number of triangles, relative to the contained number of edges. 
Implicitly, said graphs must also contain some vertices with substantially smaller degree, which is what will serve as the basis upon which the proofs are founded.
Specifically, we will use the somewhat artificial notion of \textit{edge-degrees}, formally defined as follows.
\begin{definition}
Let $G=(V,E)$ be a graph and let $e=(u,v)\in E$ be an edge. If the vertex-degrees of $u$ and $v$ are respectively $d_u$ and $d_v$, then the \textit{edge-degree} $\delta(e)$ of $e$ is $d_u+d_v$.
\end{definition}

The first infinite family of signed digraphs whose spectral characterization we will discuss is based on maximally dense clique expansions of $C_4$.
Using the results from Section \ref{sec: counting substructures}, we may be certain that any graph $G$ that is underlying to a signed digraph $D$ that is cospectral to such a clique expansion must contain precisely $\binom{n-1}{2}+1$ edges and at least $\binom{n-1}{3}-n+3$ triangles. 
Now, in case the minimum degree of $G$ is small enough, we can easily use just the number of edges to pin-point its structure precisely, which is formalized in the following lemma.
\begin{lemma}\label{lemma: low degree}
Let $n\geq 3$, and let $G$ be an $O_3$-free graph with $n$ vertices and $m=\binom{n-1}{2}+1$ edges. 
Futher, let $u$ be the vertex of minimal degree. 
If $d_u=1$, then $G\cong CE(P_3,[n-2~1~1])$. 
Moreover, if $d_u=2$ then either $G\cong CE(C_4,[n-3~1~1~1])$ or $G\cong CE(Gem,[1~n-4~1~1~1]).$
\end{lemma}
If, instead, the minimum degree exceeds two, we may instead show that the graph must be one of four exceptional graphs on small $n$. 
\begin{lemma}\label{lemma: list of c4 sufficient triangles}
Let $n\geq 3$, and let $G$ be an $O_3$-free graph with $n$ vertices, $m=\binom{n-1}{2}+1$ edges and $t = |T(G)|$ triangles. 
Then $t\geq \binom{n-1}{3}-n+3$ if and only if $G$ has minimum vertex degree $1$ or $2$, or if $G$ is one of the exceptional graphs $G_1,G_2,G_3,G_4$, illustrated in Figure \ref{fig: g1 g2 g3 g4 illustrations}. 
\end{lemma}
 \begin{figure}[t]
    \centering
     \begin{subfigure}[b]{0.32\textwidth}\centering
     \begin{tikzpicture}[scale=1]
        
        \node[group] (p) at (0,0) {$n-2$};
        \node[vertex] (1) at (1,0) {};
        \node[vertex] (2) at (2,0) {};
        
        \draw[edge] (1) to node{} (p);
        \draw[edge] (1) to node{} (2);
        
     \end{tikzpicture}
     \caption{ $CE(P_3,[n-2~1~1])$}
     
     \end{subfigure}
     \begin{subfigure}[b]{0.32\textwidth}\centering
     \begin{tikzpicture}[scale=1]
        \node[group] (1) at (1,1) {$n-3$};
        \node[vertex] (2) at (1,0) {};
        \node[vertex] (3) at (0,0) {};
        \node[vertex] (4) at (0,1) {};
        
        \draw[edge] (1) to node{} (2);
        \draw[edge] (3) to node{} (2);
        \draw[edge] (3) to node{} (4);
        \draw[edge] (1) to node{} (4);
        
     \end{tikzpicture}
     \caption{$CE(C_4,[n-3~1~1~1])$ }
     
     \end{subfigure}
     \begin{subfigure}[b]{0.32\textwidth}\centering
     \begin{tikzpicture}[scale=1]
        \node[group] (1) at (0,1) {$n-4$};
        \node[vertex] (2) at (1,1.2) {};
        \node[vertex] (3) at (0,0) {};
        \node[vertex] (4) at (1,-.2) {};
        \node[vertex] (5) at (2,.5) {};
        
        \draw[edge] (1) to node{} (2);
        \draw[edge] (1) to node{} (3);
        \draw[edge] (3) to node{} (4);
        
        \draw[edge, bend right=7] (1) to node{} (5);
        \draw[edge] (2) to node{} (5);
        \draw[edge, bend left=7] (3) to node{} (5);
        \draw[edge] (4) to node{} (5);
     \end{tikzpicture}
     \caption{$CE(Gem,[1~n-4~1~1~1])$ }
     
     \end{subfigure}
    \caption{The graphs in Lemma \ref{lemma: low degree}.}
    \label{}
\end{figure}
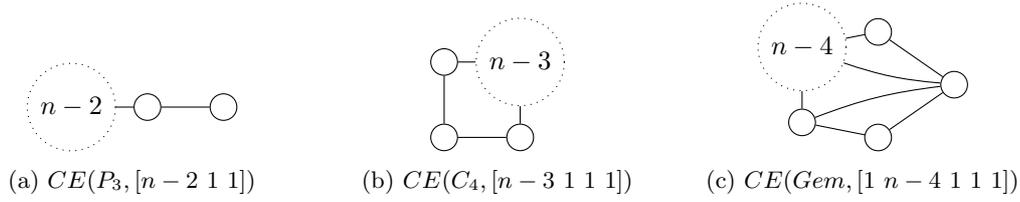
 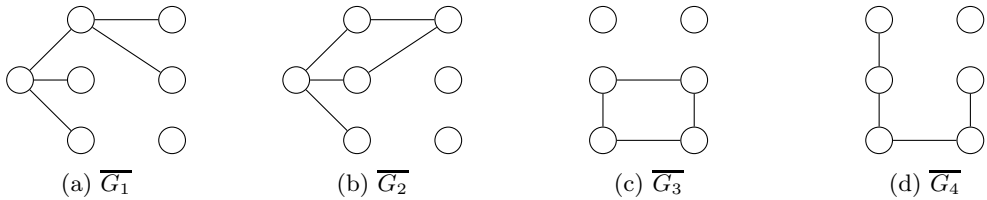
\begin{figure}[b]
    \centering
    \begin{subfigure}[b]{0.24\textwidth}\centering
     \begin{tikzpicture}[scale=.8]
        
        \node[vertex] (p) at (0,0) {};
        \node[vertex] (1) at (1,1) {};
        \node[vertex] (2) at (1,0) {};
        \node[vertex] (3) at (1,-1) {};
        \node[vertex] (4) at (2.5,1) {};
        \node[vertex] (5) at (2.5,0) {};
        \node[vertex] (6) at (2.5,-1) {};

        \draw[edge] (p) to node{} (1);
        \draw[edge] (p) to node{} (2);
        \draw[edge] (p) to node{} (3);
        \draw[edge] (1) to node{} (4);
        \draw[edge] (1) to node{} (5);
     \end{tikzpicture}
     \caption{\label{G2}$\overline{G_1}$ }
     
     \end{subfigure}
    \begin{subfigure}[b]{0.24\textwidth}\centering
    \begin{tikzpicture}[scale=.8]
        \node[vertex] (p) at (0,0) {};
        \node[vertex] (1) at (1,1) {};
        \node[vertex] (2) at (1,0) {};
        \node[vertex] (3) at (1,-1) {};
        \node[vertex] (4) at (2.5,1) {};
        \node[vertex] (5) at (2.5,0) {};
        \node[vertex] (6) at (2.5,-1) {};
        
        \draw[edge] (p) to node{} (1);
        \draw[edge] (p) to node{} (2);
        \draw[edge] (p) to node{} (3);
        \draw[edge] (4) to node{} (1);
        \draw[edge] (4) to node{} (2);
    \end{tikzpicture}
     \caption{\label{G1}$\overline{G_2}$ }
     \end{subfigure}
     \begin{subfigure}[b]{0.24\textwidth}\centering
     \begin{tikzpicture}[scale=.8]
        \node[vertex] (1) at (1,1) {};
        \node[vertex] (2) at (1,0) {};
        \node[vertex] (3) at (1,-1) {};
        \node[vertex] (4) at (2.5,1) {};
        \node[vertex] (5) at (2.5,0) {};
        \node[vertex] (6) at (2.5,-1) {};
        
        \draw[edge] (3) to node{} (2);
        \draw[edge] (5) to node{} (2);
        \draw[edge] (5) to node{} (6);
        \draw[edge] (3) to node{} (6);
     \end{tikzpicture}
     \caption{\label{G3}$\overline{G_3}$ }
     
     \end{subfigure}
     \begin{subfigure}[b]{0.24\textwidth}\centering
     \begin{tikzpicture}[scale=.8]
        \node[vertex] (1) at (1,1) {};
        \node[vertex] (2) at (1,0) {};
        \node[vertex] (3) at (1,-1) {};
        \node[vertex] (4) at (2.5,1) {};
        \node[vertex] (5) at (2.5,0) {};
        \node[vertex] (6) at (2.5,-1) {};
        
        \draw[edge] (1) to node{} (2);
        \draw[edge] (2) to node{} (3);
        \draw[edge] (3) to node{} (6);
        \draw[edge] (6) to node{} (5);
     \end{tikzpicture}
     \caption{\label{G4}$\overline{G_4}$ }
     
     \end{subfigure}
    \caption{Complements of the four exceptional graphs for Lemma \ref{lemma: list of c4 sufficient triangles}.}
    \label{fig: g1 g2 g3 g4 illustrations}
\end{figure}

\begin{proof}
First, we note that $G$ cannot have an isolated vertex, because then the remaining $n-1$ vertices would have to harbour $\binom{n-1}{2}+1$ edges, which is impossible for a simple graph. 
Thus, suppose that every vertex has degree at least $3$. 

We will then consider the triangle-free complement $\bG=(V,\bE)$ of $G$, whose degrees are denoted $d_u,~u\in V$. 
This complement has $n-2$ edges, so $\sum_{u\in V} d_u = 2n-4$. 
By the above assumption, $d_u\leq n-4$ for $u\in V$. 

By inclusion-exclusion, we may express the number of triangles in $G$ as 
\begin{equation}\label{eq: t c4}t = \binom{n}{3}- (n-2)^2 + \sum_{u\in V}\binom{d_u}{2}.\end{equation}
Here, the second term represents the $n-2$ edges missing from $G$, which are each responsible for $n-2$ missing triangles, unless they intersect with another missing edge. 
Then, using that $\bG$ is $K_3$-free, the third term corrects the overshoot resulting from the second term.

Expanding the third term and plugging in $\sum_{u\in V} d_u$ and $\sum_{u\in V}d_u^2=\sum_{e\in\bE}\delta(e)$ yields
\[t=\binom{n}{3}- (n-1)(n-2)+\frac{1}{2} \sum_{e\in\bE}\delta(e),\]
which may in turn be combined with $t\geq \binom{n-1}{3}-n+3$ to find that 
\[\sum_{e\in\bE}\delta(e)\geq n^2-5n+8.\]
Now, we may take the average of $\delta(e)$ over $\bE$ to obtain 
\begin{equation}\label{eq: delta(e) n-2}\frac{1}{n-2}\sum_{e\in\bE} \delta(e) \geq n-3 +\frac{2}{n-2} > n-3.\end{equation}
So there is an edge $e$ with $\delta(e)\geq n-2$. 
But $\delta(e)\leq n-1$, because $\bG$ contains only $n-2$ edges. 
We distinguish two cases: either there is an edge with degree $n-1$, or there is not. 

First, suppose that there is an edge $e^*=(u,v)$ with $\delta(e^*)=n-1.$ 
        Then $d_u=n-1-d_v$, and without loss of generality $d_v\leq d_u$. 
        Now, since $\bG$ is triangle-free and $u$ and $v$ are together adjacent to $n-3$ of the remaining vertices, the degree sequence of $\bG$ is $n-1-d_v, d_v, 1, 1,\ldots,1,0.$
        Observe that since $d_u\leq n-4$, it follows that $3\leq d_v\leq \frac{1}{2}(n-1)$, which we may use to we obtain an upper bound for $\sum_{e\in\bE}\delta(e)$ as:
        \[\sum_{e\in\bE}\delta(e) = \sum_{w\in V}d_w^2 = d_v^2 + (n-d_v-1)^2 + (n-3) \leq 9 + (n-4)^2 + n-3 = n^2 -7n + 22.\]
        Together with the lower bound $\sum_{e\in\bE}\delta(e)\geq n^2 - 5n + 8$, this yields $n\leq 7$, which in turn implies $d_v=d_u=3, n=7$.
        This uniquely characterizes the graph $G_1$ since $\bG$ is triangle-free and $\delta(e^*)=|E|+1.$  

Now, instead suppose that all edges have $\delta(e)\leq n-2$. 
        As argued in \eqref{eq: delta(e) n-2}, there is an edge $e^*=(u,v)$ with $\delta(e^*)=n-2$. 
        This fixes all but one edge. 
        Again, we have $d_u=n-2-d_v$, with $d_v\leq n-2-d_v\leq n-4$, so $2\leq d_v\leq \frac{1}{2}n-1$. 
        Working analogously to before, this case surprisingly yields the same upper bound for $\sum_{e\in\bE}\delta(e)$:  \[\sum_{e\in\bE}\delta(e)\leq d_v^2 + (n-2-d_v)^2 + 2\cdot 4 + n-6 \leq 4 + (n-4)^2 + n + 2 = n^2 -7n + 22. \]
        Because of the lower bound for $\sum_{e\in\bE}\delta(e),$ we then find that $6\leq n\leq 7$, and because $d_v\leq \frac{1}{2}n-1$ that $d_v=2$ and $d_u=n-4$. By checking the remaining 7 (non-isomorphic) configurations, if follows that $G$ has sufficient triangles only when it is one of the remaining exceptional graphs $G_2$, $G_3$, or $G_4$.
\end{proof}

Given the above, we may now simply consider a limited number of potential underlying graphs and investigate in which cases they lead to cospectrality.
Let $C_4^*$ be the four-cycle with gain $-\w$.
Since it belongs to the only switching isomorphism class on $C_4$ that has eigenvalue $-1$ with multiplicity $n-3$ , we may apply interlacing to conclude the following with relative ease. 
\begin{theorem}
\label{thm: c4 expansions DES}
Let $\Phi = CE(C_4^*,[n-3~1~1~1])$. Then $\Phi$ is DES, for $n\geq 4$.
\end{theorem}
\begin{proof} 
Suppose that $D$ is cospectral to $\Phi$. 
Then the spectrum of $D$ contains $2$ strictly positive and $n-2$ strictly negative eigenvalues; specifically, $-1$ occurs with multiplicity $n-3$. 
Furthermore, $|E(D)|=\binom{n-1}{2}+1$ and $|T(\Gamma(D))|\geq \binom{n-1}{3}-n+3$, which by Lemma \ref{lemma: list of c4 sufficient triangles} implies that $\Gamma(D)$ contains a vertex of degree $1$ or $2$, or is one of four exceptional graphs. 
We first explore the four exceptions, illustrated in Figure \ref{fig: g1 g2 g3 g4 illustrations}.

Recall that graphs $G_1,G_2,$ and $G_4$ contain exactly $|T(\Phi)|$ triangles, and thus $D = (G_j,\varphi)$, $j=1,2,4$, may be cospectral to $\Phi$ only when all of its triangles (i.e., a basis of its cycle space) have gain $1$. 
Using Proposition \ref{prop: we can always make a tree equal}, it follows that $D$ is switching isomorphic to its underlying graph. 
Then, simply computing the spectra of $G_1,G_2,G_4$ leads to the desired conclusion. 
If $D=(G_3,\varphi)$, then we find analogously that $6$ triangles in $D$ must have gain $1$, and exactly two must have gain $\omega$, which again leads to a unique switching isomorphism class on $G_3$, whose spectrum does not coincide with $\Phi$.
Thus, the exceptional cases are covered. 

We move on to the general case: suppose $\Gamma(D)$ contains a vertex of degree at most $2$. 
Then $\Gamma(D)$ is either $CE(P_3,[n-2~1~1])$, $CE(Gem,[1~n-4~1~1~1])$ or $CE(C_4,[n-3~1~1~1])$, by Lemma \ref{lemma: low degree}. 
Now, it may easily be brute-forced that signed digraphs 
on $Gem$ have at least $4$ eigenvalues that are not $-1$. 
By an application of eigenvalue interlacing 
, it follows that signed digraphs on clique expansions of $Gem$ have eigenvalue $-1$ with multiplicity at most $n-4$, which is insufficient. 
Similarly, signed digraphs on clique expansions of $P_3$ have eigenvalue $-1$ with multiplicity at most $n-3$; it is not hard to see that this is attained only when such a signed digraph is switching isomorphic to its underlying graph. 
Since $\trace \E(CE(P_3,[n-2~1~1]))^3>\trace \E(\Phi)^3$, it follows that $\Gamma(D)=CE(C_4,[n-3~1~1~1]).$
Finally, note that all triangles in $D$ then must have gain $1$, to satisfy $\trace \E(D)^3=\trace \E(\Phi)^3$, and all of the $4$-cycles must have gain $-\w$, since $C_4^*$ is the only 4-cycle with an eigenvalue $-1$.
Since $D$ coincides with $\Phi$ on all cycle gains, the conclusion follows by Proposition \ref{prop reff}. 
\end{proof}

\begin{table}[b]\centering
\begin{minipage}{0.40\textwidth}\centering
    \begin{minipage}{.49\textwidth}\centering
     \begin{tikzpicture} [scale = .8]
     \node[group] (1) at (-1,2){~$K_3$~~};
        \node[group] (2) at (-1,0){~$K_4$~~};
        \node[vertex] (3) at (0,1){};
        \node[vertex] (4) at (1,1){};
        
        \draw[edge] (1) to node{} (2);
        \draw[arc] (1) to node{} (3);
        \draw[edge] (2) to node{} (3);
        \draw[edge] (3) to node{} (4);
     
     \end{tikzpicture}\\
     (a)
     \end{minipage}
     \begin{minipage}{.49\textwidth}\centering
     \begin{tikzpicture} [scale = .8]
      \node[group] (1) at (0,0){~$K_6$~~};
        \node[vertex] (2) at (1,1){};
        \node[vertex] (3) at (0,2){};
        \node[vertex] (4) at (-1,1){};

        \draw[arc] (1) to node{} (2);
        \draw[edge] (2) to node{} (3);
        \draw[edge] (3) to node{} (4);
        \draw[edge] (4) to node{} (1);
     
     \end{tikzpicture}\\ \vskip3mm
     (b)
     \end{minipage}

     \captionof{figure}{A cospectral pair}
     \label{fig: c4 expansion cospecral pair}
\end{minipage}~~~~
\begin{minipage}{0.59\textwidth}\centering
    \begin{minipage}{.49\textwidth}\centering
     \begin{tikzpicture} [scale = .8]

        \node[vertex] (1) at (0,.5){};
        
        \node[vertex] (2) at (1,0){};
        \node[vertex] (6) at (1,1){};
        
        \node[vertex] (3) at (2,0){};
        \node[vertex] (5) at (2,1){};
        
        \node[vertex] (4) at (3,.5){};

        \draw[edge] (1) to node{} (2);
        \draw[edge] (1) to node{} (6);
        \draw[edge] (2) to node{} (6);
        \draw[edge] (2) to node{} (3);
        \draw[edge] (3) to node{} (4);
        \draw[edge] (3) to node{} (5);
        \draw[edge] (4) to node{} (5);
        \draw[edge] (6) to node{} (5);
     \end{tikzpicture}\\ \vskip5mm
     (a) $G_5$
     \end{minipage}
     \begin{minipage}{.49\textwidth}\centering
          \begin{tikzpicture} [scale = .8]
        
        \node[vertex] (1) at (-.7,0.5){};
        \node[vertex] (2) at (0,0){};
        \node[vertex] (3) at (1.4,0){};
        \node[vertex] (4) at (2.1,.5){};
        \node[vertex] (5) at (1.4,1){};
        \node[vertex] (6) at (.7,.5){};
        \node[vertex] (7) at (0,1){};

        \draw[edge] (1) to node{} (2);
        \draw[edge] (1) to node{} (6);
        \draw[edge] (1) to node{} (7);
        \draw[edge] (2) to node{} (3);
        \draw[edge] (2) to node{} (6);
        \draw[edge] (2) to node{} (7);
        \draw[edge] (3) to node{} (4);
        \draw[edge] (3) to node{} (5);
        \draw[edge] (3) to node{} (6);
        \draw[edge] (4) to node{} (5);
        \draw[edge] (5) to node{} (6);
        \draw[edge] (6) to node{} (7);
        \end{tikzpicture}\\  \vskip5mm
     (b) $G_6$
     \end{minipage}

     \captionof{figure}{Exceptional graphs for Lemma \ref{lemma: list of graphs with sufficient triangles }}
     \label{fig: rewrite exceptions}
\end{minipage}
\end{table}

Moreover, note that as a consequence of the proof above, we get the following for free, since all  candidates for cospectrality that contain sufficient triangles are switching isomorphic.
\begin{proposition}\label{prop: p3 expansion DES}
$CE(P_3,[n-2~1~1])$ is DES for $n\geq 3$.
\end{proposition}

At this point, the attentive reader may wonder whether the above holds analogously for clique expansions of the other switching classes on $C_4$. 
While most of the argument will hold up, we find that there are many signed digraphs on expansions of $P_3$ that have an eigenvalue $-1$ with sufficiently high multiplicity to potentially share the spectrum of such a $C_4$-expansion. 
In fact, an example of such a cospectral pair is shown in Figure \ref{fig: c4 expansion cospecral pair}. 
Thus, we would have to provide a substantially different approach; in the interest of unity, we move on to the next family of graphs. \\

In similar fashion, we may also consider maximally dense expansions of $C_5$, that satisfy $\lambda_2>0>\lambda_3$. 
We follow largely the same line as in the proof of Lemma \ref{lemma: list of c4 sufficient triangles}. 
\begin{lemma}\label{lemma: list of graphs with sufficient triangles }
Let $G$ be an $O_3$-free graph with $n\geq 5$ vertices, $m=\binom{n-2}{2}+2$ edges and $t=|T(G)|$ triangles. Then $t\geq \binom{n-2}{3}-n+4$ if and only if $G$ is one of the following graphs:
\begin{enumerate}[label=\roman*)]
    \item $CE(C_5,[n-4~~1~~1~~1~~1])$,
    \item $CE(P_4,[n-3~~1~~1~~1])$, i.e., an $(n-2,2)$-kite, 
    \item $CE(P_4,[2~~1~~n-4~~1])$,
    \item One of the sporadic examples $CE(P_3,[3~~1~~3])$, $G_5$ or $G_6$. (See Figure \ref{fig: rewrite exceptions}.)
\end{enumerate}
\end{lemma}
\begin{proof}
Like before, $G$ cannot have an isolated vertex. 
If $G$ has a vertex of degree $1$, then the $n-2$ non-neighbors form a clique, and it follows that $G$ is an $(n-2,2)$-kite. 
For the remaining cases, we may assume that every vertex has degree at least $2$. 

Again consider $\bG$, whose degrees are $d_u,~u\in V$, and note that now $\sum_{u\in V} d_u = 4n-10$. 
Moreover, by assumption, $d_u\leq n-3$ for $u\in V$. 
As before, by inclusion-exclusion, it follows that the number of triangles $t$ in $G$ may be expressed as
\[t = \binom{n}{3}-(2n-5)(n-2) + \sum_{u\in V}\binom{d_u}{2}.\]
Using that $\sum_{e\in\bE} \delta(e) = \sum_{u\in V} d_u^2$ and the above, it follows that \[t\geq \binom{n-2}{2} -n+4\iff\sum_{e\in\bE} \delta(e) \geq 2n^2 - 8n + 10.\]
From this inequality, we will determine the several remaining options for $\bG$.
Once more, we take the average over $\bE$ to see that
\[\frac{1}{2n-5}\sum_{e\in\bE} \delta(e) \geq n-2 + \frac{n}{2n-5} > n-2.\]
This implies that there is an edge $e^*\in \bE$, such that $\delta(e^*)\geq n-1.$
Note however that, because $\bG$ is triangle-free, $\delta(e)\leq n$.
Now, either (I) there is an edge with degree $n$ or (II) there is not.

(I): First, assume that $\bG$ has an edge $e^* = (u,v) \in \bE$ such that $\delta(e^*)=n$ and $d_u=n-d_v$.
Without loss of generality, $d_v\leq d_u$, and because $n-d_v\leq n-3$, it follows that $3\leq d_v\leq \frac{1}{2}n$.

Let $\bE'$ be the set of $n-1$ edges in $\bE$ that are incident to $u$ or $v$, let $N_u = N(u)\setminus \{u,v\}$ be the set of neighbors of $u$ besides $v$ and similarly $N_v=N(v)\setminus \{u,v\}.$
Because there are $n-4$ edges not in $\bE'$, and because $\bG$ is triangle-free, these are edges with one vertex in $N_u$ and the other in $N_v$. 
Thus, $\sum_{w\in N_u}( d_w-1) =\sum_{w\in N_v}(d_w-1) = n-4$. 
Now 
\begin{align*}
    \sum_{e\in \bE'}\delta(e)    &= \delta(e^*) + \sum_{w\in N_u}(d_w + n-d_v) + \sum_{w\in N_v}(d_w + d_v)\\
                            &= n^2 - (2d_v-3)n + 2d_v^2 - 10. 
\end{align*}
So for the remaining $n-4$ edges not in $\bE'$ we require that $\sum_{e\not\in\bE'} \geq n^2 + (2d_v-11)n - 2d_v^2 + 20.$

Suppose now that $d_v\geq 4$. Because the lower bound $n^2 + (2d_v-11)n - 2d_v^2 + 20$ is the weakest for $d_v=4$ (in the range $4\leq d_v\leq \frac{n}{2}$), we obtain that \[\sum_{e\not\in\bE'}\delta(e) \geq n^2 -3n - 12\text{, and thus }\frac{1}{n-4}\sum_{e\not\in\bE'}\delta(e)\geq n-1 + \frac{2n-16}{n-4}.\]
This implies that the only possible values for $n$ and $d_v$ in this range is $n=8$, $d_v=4$, in which case $\delta(e)=n-1$ for all $e\in \bE'$. 
It is easy to check however that this requires more edges than $\bE'$ contains. 
In addition, it is easy to check that if $e\not\in \bE'$ then $\delta(e)\leq n-1$, for otherwise the complement of $\bE'$ would consist of at least $n-3$ edges. 

So instead we must have $d_v=3$ and $\sum_{e\not\in \bE'} \delta(e) \geq n^2 - 5n + 2$. 
Averaging yields \[\frac{1}{n-4}\sum_{e\not\in\bE'}\delta(e) \geq n-2 + \frac{n-6}{n-4},\] which implies that if $n>6$ then there is an edge $\tilde{e}\not\in\bE'$ with $\delta(\tilde{e})=n-1$. 
Let us assume existence of $\tilde{e}$ for $n=6$ as well, and let $\tilde{e}=(\tilde{u},\tilde{v}),$ with $\tilde{u}\in N_v$ and $\tilde{v}\in N_u$. 
Then there are two options. 
The first is that $\delta_{\tilde{v}}=2$ and $\delta_{\tilde{u}}=n-3.$ 
This fixes all edges not in $\bE'$, and we obtain the complement of $CE(P_4,[2~~1~~n-4~~1])$. 
The second option is that $d_{\tilde{v}}=3$ and $d_{\tilde{u}} = n-4$, where we may assume that $n>6$, otherwise this is the same as the first option. 
Again, this fixes the entire graph.
However only for $n=7$ does it satisfy the requirements, and we obtain the complement of $G_6$. 
For $n=6$, the above inequality does not guarantee the existence of an edge $\tilde{e}\not\in \bE'$ with $d(\tilde{e})=n-1$. 
Indeed, if we require all edges $e\not\in\bE'$ to have $\delta(e)\leq n-2$, then we obtain the complement of $G_5$. 

(II): Let us consider the possibility that, contrary to the case above, there is no edge $e\in\bE$ with $\delta(e) = n,$ but there is an edge $e^*=(u,v)$ such that $\delta(e^*)=n-1$. 
Again, let $d_u=n-1-d_v$ and assume without loss of generality that $d_v\leq d_u$. 
Since $n-1-d_v\leq n-3$, we now have that $2\leq d_v\leq \frac{1}{2}(n-1).$
We proceed in the same way as before, but now $\bE'$ has $n-2$ edges. 
The main difference to (I) is that there is a vertex $z$ that is not adjacent to $u$ or $v$. That is: $V\setminus (N_u\cup N_v) = \{u,v,z\}.$
In this case, we have that $\sum_{w\in N_u} (d_w-1) + \sum_{w\in N_v}(d_w-1) + d_z = 2(n-3)$ and therefore
\begin{align*}
    \sum_{e\in\bE'}\delta(e) &= n-1 + (n-d_v-2)(n-d_v-1) + (d_v-1)d_v + \sum_{w\in N_u} d_w + \sum_{w\in N_v}d_w \\ &= n^2 -(2d_v-1)n + 2d_v^2 + 2d_v - 8-d_z. 
\end{align*}
It follows that for the remaining $n-3$ edges, we require 
\[\sum_{e\not\in\bE'} \delta(e) \geq n^2 + (2d_v-9)n - 2d_v^2 -2d_v + 18 + d_z.\]
By averaging (as before) over the $n-3$ edges and using that this average is at most $n-1$, we find that we may restrict to the cases $d_v=3$ with $7\leq n\leq 9$ and $d_v=2$ with $n\geq 5$. 
The latter case leads (as only possibility) to the complement of $CE(C_5,[n-4~~1~~1~~1~~1]).$

For $d_v=3$, we have $\sum_{e\not\in \bE'} \delta(e) \geq n^2 -3n - 6 + d_z$ and $d_w \leq n-4$ for all $w\in V\setminus \{z\}$.
If $d_z\geq 1$, then consider an edge $e'$ incident to $z$.
It must satisfy $\delta(e')\leq d_z + n-4$, and thus \[\sum_{e\not\in E'\cup\{e'\}}\delta(e)\geq n^2 - 4n - 2.\]
But then $\frac{1}{n-4}\sum_{e\not\in E'\cup\{e'\}}\delta(e)\geq n-1,$ which is a contradiction.

Therefore, $d_z=0$, and we have to add $n-3$ edges between $N_u$ and $N_v$. 
For $n=7$, this yields $K_1 \cup K_{3,3}$, which is the complement of our final sporadic example $CE(P_3,[3~1~3])$.
For $n=8,9$, it is easily verified that the corresponding graphs (respectively the complements of $K_1\cup K_{3,4}$ minus one edge, and $K_1\cup K_{3,5}$ minus two edges) contain insufficient triangles, which concludes the proof.
\end{proof}

\begin{theorem}\label{thm: c5 DES}
Let $G=CE(C_5,\tau)$ with $\tau = [n-4~~1~~1~~1~~1]$, and let $\Phi=(G,\varphi)$  where $\varphi$ is of type $A$ or $C$. Then $\Phi$ is DES. 
\end{theorem}
\begin{proof}
Suppose that $D$ is cospectral to $\Phi$.
Using Lemma \ref{lemma: triangles}, it follows that $\Gamma(D)$ contains at least $\frac{1}{6}\trace \E(\Phi)^3=\binom{n-2}{3}-n+4$ triangles, which by Lemma \ref{lemma: list of graphs with sufficient triangles } implies that $\Gamma(D)$ is one of at most six potential graphs. 
Analogously to the proof of Theorem \ref{thm: c4 expansions DES}, it may easily be evaluated that the sporadic cases are not underlying to any signed digraphs that are cospectral to appropriately sized $\Phi$, so we focus on the general case. 

Suppose that $\Gamma(D)=CE(P_4,[n-3~~1~~1~~1])$; an $(n-2,2)$-kite.
By Corollary \ref{cor: kite graphs}, the signed digraph that is induced by the $(n-2)$-clique is switching isomorphic to $K_{n-2}$ or $K_{n-2}^*$.
This implies that $\frac{1}{6}\trace \E(D)^3 \geq \binom{n-2}{3} - \frac{1}{2}(n-4)>\frac{1}{6}\trace \E(\Phi)^3,$ which is a contradiction. 

Next, suppose that $\Gamma(D)=CE(P_4,[2~~1~~n-4~~1])$. 
Then $\Gamma(D)$ contains exactly one more triangle than $G$ for every $n\geq 5$, which implies that $\trace \E(D)^3 \not= \trace \E(\Phi)^3$, and we again obtain a contradiction. 
Hence, $\Gamma(D)=CE(C_5,[n-4~~1~~1~~1~~1])$, and the conclusion regarding $\varphi$ follows by Proposition \ref{prop reff}. 
\end{proof}
As before, we can draw similar conclusions for other graphs considered above.
\begin{theorem}\label{thm: kite DES}
Let $G$ be an $(n-2,2)$-kite, $n\geq 3$, and let $\Phi=(G,\varphi)$ be such that the induced $(n-2)$-clique is switching isomorphic to either $K_{n-2}$ or $K_{n-2}^*$. 
Then $\Phi$ is DES.
\end{theorem}
\begin{proof}
By Corollary \ref{cor: kite graphs}, $\Phi$ has $\lambda_2>0>\lambda_3$, so Lemma \ref{lemma: list of graphs with sufficient triangles } is applicable if $n\geq 5$.
Now, observe that for $n>6$, $\frac{1}{6}\trace \E(\Phi)>\binom{n-2}{3}-n+5$, which is the largest number of triangles in any graph in Lemma \ref{lemma: list of graphs with sufficient triangles } that is not itself an $(n-2,2)$-kite. 
The conclusion follows easily by brute-forcing (by computer) the limited collection of signed digraphs from Lemma \ref{lemma: list of graphs with sufficient triangles } on $n=4,5,6$ whose underlying graphs do contain sufficient triangles. (Recall that if $n=3$ then $\Phi\sim G$, by Proposition \ref{prop: we can always make a tree equal}.)
\end{proof}
Finally, as an immediate consequence of Theorems \ref{thm: c5 DES} and \ref{thm: kite DES}, the following is easily verified.
Since all cospectral candidates for $n\geq 8$ are DES themselves, one only needs to check a limited number of graphs on $5\leq n\leq 7.$ 
\begin{proposition}
$CE(P_4,[2~1~n-4~1])$ is DES.
\end{proposition}
To conclude, we note that there are various families of signed digraphs that are close tangents of the discussed DES families, which have remained untreated in this section. 
For example, maximally dense $C_5$ expansions with signatures $B$ and $D$, or minimally connected semi-complete graphs come to mind. 
While the authors are convinced that similar results could be obtained for these cases, their particular challenges are preserved for future research. 

\section{Open questions}
\label{sec: open}
We end this article with a summarizing list of open questions. 
\begin{question}
For a given underlying graph, there is a switching equivalence class whose edges have gain $1$ on a predetermined spanning tree. However, is it possible to determine (or bound) the number of members of said class that do so?
\end{question}
\begin{question}
An interesting, 'ugly' example of a signed digraph whose spectrum is symmetric was provided in Figure \ref{fig: ugly example}.
Perhaps, it is too ambitious to ask for a tight characterization of all signed digraphs with symmetric spectra. However, would it be possible to formulate some necessary properties?
\end{question}
\begin{question}
The collections of signed digraphs whose rank is either $2$ or $3$ have been completely characterized.
Do ranks $4$ and up also allow for a comprehensive characterization?
\end{question}
\begin{question}
In order to keep the discussion on signed digraphs with two non-negative eigenvalues tidy, we have zoomed in on a number of special cases. 
Is it feasible to complete the characterization? Possibly under some additional restrictions?
\end{question}
\begin{question}
How do the results of this work change when one considers other gain groups $\mathbb{T}_k$? Which parts carry over to a gain set that is not closed under multiplication? 
\end{question}
\begin{question}
Do the remaining columns of Table \ref{tab: C5 lex expac vectors} give rise to DES signed digraphs? And what about other signed digraphs on $C_4$?
\end{question}
\begin{question}
In our representation, $\varphi(u,v)$ is positive if and only if $\varphi(v,u)$ is positive. 
In dynamic systems, this does not have to be the case. 
Should we adapt our formulation to accommodate for these 'negative feedback loops', or can the issue be sufficiently circumvented with a simple subdivision trick? 
Moreover, can spectral analysis of $\mathbb{T}_6$-gain graphs then offer new insights into the (possible) stability of such systems?
\end{question}

\section*{Acknowledgements}
The authors would like to express their gratitude to Jack Koolen and Willem Haemers for their invaluable comments in weighing the advantages and disadvantages of various  Hermitian adjacency matrices for directed graphs. 
Furthermore, the anonymous referees are thanked for their valuable remarks and suggestions to the betterment of this work.

\scriptsize
\DeclareRobustCommand{\VAN}[3]{#3}
\bibliographystyle{abbrvnat}
\bibliography{mybib}{}
\normalsize
\end{document}